
\magnification=\magstephalf
\hoffset=0.5truein
\input amstex
\documentstyle{amsppt}
\refstyle{A}
\loadeusm

\def\cI{{\Cal I}}

\def\bbR{{\Bbb R}}\def\bbC{{\Bbb C}}\def\bbB{{\Bbb B}}
\def\bbZ{{\Bbb Z}}\def\bbP{{\Bbb P}}

\def\w{{\mathchoice{\,{\scriptstyle\wedge}\,}{{\scriptstyle\wedge}}
      {{\scriptscriptstyle\wedge}}{{\scriptscriptstyle\wedge}}}}
\def\eusu{\operatorname{\frak su}}
\def\eusl{\operatorname{\frak sl}}
\def\eug{\operatorname{\frak g}}
\def\SL{\operatorname{SL}}
\def\SU{\operatorname{SU}}\def\Un{\operatorname{U}}

\def\Lie{\operatorname{\eusm L}}
\def\ts{\textstyle }
\def\>{\medspace}
\def\la{\langle}\def\ra{\rangle}

\def\gbold{{\bold g}}
\def\ebold{{\bold e}}
\def\vbold{{\bold v}}

\def\abold{{\bold a}}
\def\sbold{{\bold s}}
\def\tbold{{\bold t}}
\def\pbold{{\bold p}}
\def\qbold{{\bold q}}
\def\ibold{{\bold i}}
\def\jbold{{\bold j}}
\def\kbold{{\bold k}}

\def\bone{{\bold 1}}

\def\Abold{{\bold A}}
\def\Bbold{{\bold B}}
\def\Fbold{{\bold F}}

\def\sR{\eusm R}

\def\p{\operatorname{\frak p}}  

\topmatter
\title Levi-flat Minimal Hypersurfaces\\
in Two-dimensional Complex Space Forms
\endtitle
\author
Robert L. Bryant
\endauthor
\affil
Duke University
\endaffil
\address
Department of Mathematics,
Duke University,
Box 90320,
Durham, NC 27708-0320
\endaddress
\date
September 25, 1999
\enddate

\email
bryant\@math.duke.edu
\endemail

\keywords 
 Levi form, minimal, complex space forms, CR hypersurfaces
\endkeywords

\subjclass 
 Primary:   32F25. 
 Secondary: 53C42  
\endsubjclass

\dedicatory
\enddedicatory

\thanks
This article was begun during a June 1999 visit at the IHES.  I would like 
to thank the IHES for its hospitality.  The idea for this article came
to me during a conversation with Mikhail~Gromov and Gennadi~Henkin, 
who asked (perhaps idly) whether any nontrivial examples of the kind 
mentioned in the title exist. I thank them for their stimulating 
conversation.
\endthanks

\abstract
The purpose of this article is to classify the real hypersurfaces 
in complex space forms of dimension~$2$ that are both Levi-flat and
minimal.  The main results are as follows: 

When the curvature of the complex space
form is nonzero, there is a 1-parameter family of such hypersurfaces.  
Specifically, for each one-parameter subgroup of the isometry group 
of the complex space form, there is an essentially unique example that 
is invariant under this one-parameter subgroup. 

On the other hand, when the curvature of the space form is zero, i.e.,
when the space form is~$\bbC^2$ with its standard metric, there is an 
additional `exceptional' example that has no continuous symmetries
but is invariant under a lattice of translations.  Up to isometry
and homothety, this is the unique example with no continuous 
symmetries.  
\endabstract

\rightheadtext{Levi-flat minimal hypersurfaces}

\endtopmatter

\document

\head 0. Introduction  \endhead

A real hypersurface~$\Sigma^3\subset\bbC^2$ is {\it Levi-flat\/}~\cite{CM} 
if it is foliated by complex curves.  (If such a foliation exists,
it is necessarily unique.)  Thus, a Levi-flat hypersurface
in~$\bbC^2$ is essentially a 1-parameter family of complex curves 
in~$\bbC^2$.  If one imposes the further
condition that the hypersurface be {\it minimal\/}, there is,
in addition to the obvious example of a real hyperplane,
the deleted cone~$C^*\subset\bbC^2\setminus\{(0,0)\}$ defined by
$$
|z_1|^2 - |z_2|^2 = 0,\qquad\qquad |z_1|^2 + |z_2|^2 > 0.
$$
This cone is foliated by the (punctured) lines~$z_1 = \lambda\,z_2$ 
with~$|\lambda| =1$ and hence is Levi-flat. Since~$C^*$ is the cone 
on the Clifford torus, it is also minimal as a submanifold of~$\bbC^2$.

It is not obvious that there are any examples of minimal, Levi-flat
hypersurfaces in~$\bbC^2$ that are distinct from these up to rigid motion.  
The condition of being either minimal or Levi-flat constitutes
a single non-linear second order PDE for the hypersurface~$\Sigma$.
A short calculation shows that the combined conditions form a second
order system that is not involutive in Cartan's sense.  In fact,
by Cartan's classification~\cite{Ca} of the involutive second order 
systems for one function of three variables, there is no second order 
equation that is in involution with the minimal hypersurface equation
for a hypersurface in a Riemannian 4-manifold.  Thus, describing 
the solutions of such a system requires analysis
that goes beyond an application of the Cartan-K\"ahler theorem.

In this article, I carry out this analysis, classify the solutions
of this overdetermined system, both locally and globally,
and show that there are many other examples.  Since it is no 
harder to do the analysis for the general two-dimensional
complex space form, I do the computations in this more 
general setting.  While the calculations were guided by 
certain concepts from exterior differential systems, 
this article has been written so that no knowledge of this subject is
required of the reader beyond the (elementary) Frobenius theorem
on integrable plane fields.  Nevertheless, the reader who wonders
how some of the calculations in~\S2 could be motivated might want 
to consult~\cite{BCG, Chapter~VI}.  General references on 
calculations via the moving frame could also be helpful, in which case
I recommend~\cite{Sp} or~\cite{Gr}. 

The results can be described as follows:  
Each local solution extends to a unique maximal solution and  
the space of maximal solutions is finite dimensional, breaking 
up into two or three different families.  

The members of the first family are those hypersurfaces~$\Sigma$
whose complex leaves are totally geodesic in the ambient space form.  
In flat space, there are only two such examples up to isometry: 
the hyperplane and the Clifford cone constructed above.  
When the space form has positive sectional curvature
and hence is~$\bbP^2$ with its standard Fubini-Study metric up to
a constant scale factor, there is only one example up to rigid
motion.  Its closure in~$\bbP^2$ has one singular point, near which it 
resembles the Clifford cone in flat space.  When the space form has negative 
sectional curvature and hence is the complex hyperbolic $2$-ball~$\bbB^2$ 
(i.e., the noncompact dual of~$\bbP^2$) up to a scale factor, there are 
three distinct examples up to isometry.  The closure of one of these examples 
has a singular point, near which it resembles the Clifford cone.  
The other two examples are nonsingular, complete, embedded hypersurfaces. 
For details, see~\S3.1.

The remaining two families are somewhat more difficult to describe explicitly.  
The structure equations for the second family show that each
such example~$\Sigma^3$ is invariant under a one-parameter group of isometries
of the ambient space and that this one-parameter group acts on the
hypersurface~$\Sigma$ preserving each of its complex leaves.  
Conversely, each one-parameter group of isometries of the ambient space
preserves a family of holomorphic curves that foliates the ambient space
in the complement of the fixed point set.  Up to ambient isometry, 
there is a unique one-parameter family of these curves whose union is a 
minimal hypersurface.  The minimal Levi-flat hypersurfaces constructed
in this way that do not belong to the first family constitute the 
members of the second family.  In~\S3.2, I construct these hypersurfaces
explicitly for each conjugacy class of one-parameter subgroup of the
isometry group of the ambient space form.
The examples in this second family often have some sort of singular locus 
and can be either real algebraic or transcendental, see~\S3.2.

The third family is the most difficult to describe explicitly.  It only
exists when the ambient curvature is zero, i.e., in the case of~$\bbC^2$
itself.  Up to holomorphic isometry and homothety, there is
only one such example and it is periodic with respect to a lattice~$\Lambda
\subset\bbC^2$ of type~$F_4$.  
The quotient hypersurface~$\Sigma^3\subset\bbC^2/\Lambda$
has quite interesting properties.  Its complex leaves
are compact Riemann surfaces of genus~$3$ and the $1$@-parameter family of
genus 3 surfaces that makes up this hypersurface is a nontrivial variation 
in moduli.  The formula that defines the embedding of~$\Sigma$ 
into the abelian variety~$\bbC^2/\Lambda$ is essentially a quotient of the 
Abel-Jacobi mapping on each complex leaf.  There is reason to believe
that this hypersurface is an open dense subset of a `real algebraic' 
hypersuface in the algebraic variety~$\bbC^2/\Lambda$, but I have not 
verified this in detail.  I would like to thank Dave Morrison for a helpful
conversation about the algebraic geometry of this example.

\head 1. Two-Dimensional Complex Space Forms  \endhead

This section introduces the structure equations for 
complex space forms of dimension 2 and establishes the notation that will
be used for the remainder of the article.  For further discussion of 
these models, the reader might consult~\cite{He} or \cite{KN}.  

\subhead 1.1. The group~$G_R$ \endsubhead
Let~$R$ be a real number and let~$G_R\subset\SL(3,\bbC)$ be the
connected subgroup whose Lie algebra~$\eug_R$ consists of the matrices 
of the form
$$
\pmatrix
   ir_1&-R\,\bar x &-R\,\bar y\\
    x & ir_2 & -\bar z \\
    y & z    & -i(r_1+r_2)\\
\endpmatrix
$$
where~$r_1$ and $r_2$ are real and $x$, $y$, and $z$ are complex. 

When~$R\not=0$, this is the identity component of the set of unimodular
matrices~$\gbold$ that satisfy~$^t\bar\gbold H_R\gbold = H_R$, where
$$
H_R = {{}^t\overline{H_R}} = \pmatrix 1&0&0\\0&R&0\\0&0&R\\ \endpmatrix.
$$
In this case, $H_R$ defines a nondegenerate Hermitian inner 
product~$\la,\ra_R$ on~$\bbC^3$.  Even the matrix $H_0$ defines a 
(very degenerate) Hermitian inner product~$\la,\ra_0$ 
on~$\bbC^3$ and~$G_0$ preserves it.

\subhead 1.2. The complex space form~$\bbP^2_R$ \endsubhead
The set~$\bbP^2_R\subset\bbP^2$ consisting of the lines 
through~$0\in\bbC^3$ on which~$\la,\ra_R$ is positive 
is a homogeneous space of~$G_R$.  
Write the general element of~$G_R$ as 
$$
\gbold = (\ebold_0\quad \ebold_1\quad \ebold_2),
$$
where the columns~$\ebold_i$ of~$\gbold$ are to be regarded as 
$\bbC^3$-valued functions on~$G_R$.  The map~$\pi:G_R\to\bbP^2_R$ defined 
by~$\pi(\gbold) = \bbC{\cdot}\ebold_0$ is a submersion.  The fibers 
of~$\pi$ are the left cosets of the connected subgroup~$K\subset G_R$ 
whose Lie algebra consists of matrices of the form
$$
\pmatrix
   ir_1&0 &0\\
    0& ir_2 & -\bar z \\
    0 & z    & -i(r_1+r_2)\\
\endpmatrix.
$$
The group~$K$ is compact and isomorphic to the nontrivial double cover 
of~$\Un(2)$. In particular, $\bbP^2_R\simeq G_R/K$ as a homogeneous space.

\subhead 1.3. The structure equations \endsubhead
Write the left invariant Maurer-Cartan form on~$G_R$ in the form
$$
\gamma = \gbold^{-1}\,d\gbold = 
\pmatrix
   i\tau &-R\,\bar\eta &-R\,\bar\omega \\
    \eta & i(\phi{+}\tau) & -\bar \sigma \\
    \omega & \sigma    & -i(\phi{+}2\tau)\\
\endpmatrix,
$$
so that the {\it first structure equation\/} becomes
$$
(d\ebold_0\quad d\ebold_1\quad d\ebold_2) 
= (\ebold_0\quad \ebold_1\quad \ebold_2)
\pmatrix
   i\tau &-R\,\bar\eta &-R\,\bar\omega \\
    \eta & i(\phi{+}\tau) & -\bar \sigma \\
    \omega & \sigma    & -i(\phi{+}2\tau)\\
\endpmatrix\,.
$$  

There exist on~$\bbP^2_R$ a unique metric~$ds^2$
and a $ds^2$@-orthogonal complex structure~$J$ 
with corresponding K\"ahler form~$\Omega$ for which
$$
\pi^*\bigl(ds^2\bigr) = \eta\circ\bar\eta + \omega\circ\bar\omega
\qquad\text{and}\qquad
\pi^*(\Omega) = {\ts{\frac i2}}(\eta\w\bar\eta + \omega\w\bar\omega).
$$

The {\it second structure equation\/}~$d\gamma=-\gamma\w\gamma$ 
shows that this K\"ahler structure has constant holomorphic sectional 
curvature~$4R$. (I.e., the Gauss curvature 
of any totally geodesic complex curve in~$\bbP^2_R$ is~$4R$.) 

From now on, the fibration~$\pi:G_R\to\bbP^2_R$ will be taken as the 
standard unitary bundle structure for the K\"ahler geometry of~$\bbP^2_R$.  
(Strictly speaking, of course, this is not quite correct since one should 
first divide out by the center of~$G_R$, a cyclic subgroup of order~$3$, 
but for simplicity, I will not do this.  It should not cause any confusion.)

\head 2.  Real Hypersurfaces \endhead

Let~$\Sigma^3$ be a connected, smooth, embedded%
\footnote{All these calculations will be local, so embeddedness is 
not a serious restriction.}
real hypersurface in~$\bbP^2_R$.  
The  pre\-image~$B_0=\pi^{-1}(\Sigma)$ is a principal 
$K$@-bundle over~$\Sigma$. 
From now on, all the forms on~$G_R$ are to be understood as pulled back 
to~$B_0$. 

\subhead 2.1. First invariants \endsubhead
Since~$\Sigma$ is a hypersurface, there will be one linear relation among 
the real and imaginary parts of the two 1-forms~$\eta$ and~$\omega$.  
Let~$B_1\subset B_0$ be the subset where this relation is~$\eta=\bar\eta$.  
Then~$B_1$ is a union of left $K_1$@-cosets where~$K_1\simeq S^1$ is the 
group of matrices of the form
$$
E_\theta = \pmatrix e^{i\theta}&0&0\\
                    0&e^{i\theta}&\\
                    0&0&e^{-2i\theta}\\ 
           \endpmatrix.
$$
From now on, all the forms on~$B_0$ are to be understood as pulled back 
to~$B_1$.   In addition to the relation~$\eta=\bar\eta$, 
there will be relations of the form
$$
\align
  \phi &= H\,\eta - \phantom{iL}\llap{$ia$}\,\omega 
               + \phantom{2\,s}\llap{$i\bar a$}\,\bar\omega\\
\sigma &=  \phantom{H}\llap{$c$}\,\eta + iL\,\omega - 2\,s\,\bar\omega
\endalign
$$
for some functions~$a$, $c$, $H= \bar H$, $L$, and $s$ on~$B_1$.  
(The choice of numerical coefficients is cosmetic.) 
By the structure equations,
$$
d\eta = -i\,\phi\w\eta + \bar\sigma\w\omega 
= -\bigl(a\omega - \bar a\,\bar\omega)\w\eta
+\overline{\bigl(c\,\eta + iL\,\omega - 2\,s\,\bar\omega)}\w\omega.
$$
Since~$\eta$ is real, the imaginary part of the right hand expression
must vanish.   I.e.,
$$
L=\bar L \qquad\text{ and} \qquad c=-2\,\bar a.
$$  

Let $R_\theta:B_1\to B_1$ denote right action by the matrix~$E_\theta$. 
Then~$\eta$, $H$, and $L$ are invariant under~$R_\theta$ while
$$
R_\theta^*\,\omega = e^{3i\theta}\omega,\quad 
R_\theta^*\,a = e^{-3i\theta}\,a, \quad R_\theta^*\,s = e^{6i\theta}\,s .
$$
Note that quantities such as~$\eta$,  $a\,\omega$, $\bar s\,\omega^2$, $L$, 
$H$, $|a|^2$, and $|s|^2$ are $\pi$-semibasic and  invariant under~$R_\theta$ 
and so can be considered to be well defined as functions or 
forms on~$\Sigma$.

\subsubhead 2.1.1. Levi-flatness \endsubsubhead
The equation~$\eta=0$ defines the preimage in~$B_1$ of the bundle of
complex tangent spaces to~$\Sigma$.  Consequently, $\Sigma$ will be Levi-flat 
if and only if~$\eta\w d\eta=0$.  However, by the structure equations 
and the relations just derived,
$$
\eta\w d\eta =  iL\,\eta\w\omega\w\bar\omega.
$$
Thus, Levi-flatness is equivalent to the condition~$L=0$.  From now on, 
I will assume that~$\Sigma$ is Levi-flat.

\subsubhead 2.1.2. Minimality \endsubsubhead
The induced metric on~$\Sigma$ pulls 
back to~$B_1$ to be the quadratic form~$\eta^2+\omega\circ\bar\omega$, 
while the second fundamental form~$\text{I\!I}$ satisfies
$$
\pi^*(\text{I\!I})= c_1H\,\eta^2 +\text{Re}(c_2\,a\,\omega)\circ\eta
             +\text{Re}(c_3\,\bar s\,\omega^2)
$$
for some nonzero constants~$c_1$, $c_2$, and $c_3$ (the explicit values 
will not be important for what follows).  In particular, $H$ is the the mean 
curvature function of~$\Sigma$ (up to some universal constant multiple),
i.e., $\Sigma$ is minimal if and only if~$H$ vanishes identically on~$B_1$.  
From now on, I will assume that~$\Sigma$ is minimal (as well as Levi-flat).

\subhead 2.2. Differential consequences of the structure equations \endsubhead
At this point, the forms on~$B_1$ satisfy 
the reality condition~$\eta=\bar\eta$, 
the nondegeneracy condition~$\eta\w\omega\w\bar\omega\not=0$, 
and the relations
$$
\align
  \phi &= - ia\omega + i\bar a\,\bar\omega\,,\\
\sigma &=  -2\,\bar a \,\eta - 2\,s\,\bar\omega\,.
\endalign
$$
Thus, $\gamma$ pulled back to~$B_1$ has the form
$$
\gamma = 
\pmatrix
   i\tau &-R\,\eta &-R\,\bar\omega \\
    \eta & i\tau + a\,\omega - \bar a\,\bar\omega 
                       & 2\,a\,\eta + 2\,\bar s\,\omega \\
    \omega & -2\,\bar a\,\eta - 2\,s\,\bar\omega 
                & -2i\tau - a\,\omega+\bar a\,\bar\omega\\
\endpmatrix.
$$
The structure equation~$d\gamma=-\gamma\w\gamma$ 
expands to the relations
$$
\align
d\tau &= iR\,\omega\w\bar\omega,\\
d\eta &= (a\,\omega + \bar a\,\bar\omega)\w\eta,\\
d\omega &= (3i\tau - \bar a\,\bar\omega)\w\omega + 2\,s\,\bar\omega\w\eta,\\
\endalign
$$
and implies the existence of complex-valued functions~$x$ and $y$ on~$B_1$
so that
$$
\align
da &= -3ia\,\tau-6\bar a\bar s\,\eta + (x-3a^2)\,\omega
            -({\ts\frac12}R{-}|a|^2{+}2|s|^2)\,\bar\omega,\\
ds &= \phantom{-}6is\,\tau+\bar x\,\eta + 3sa\,\omega + y\,\bar\omega.\\
\endalign
$$

\remark{Remark 1}
These equations imply strong conditions about the vanishing locus of~$s$
on each complex leaf $L\subset\Sigma$.  In a small neighborhood~$U$ 
of any point~$p\in L$, one can choose a complex coordinate~$z$ so that, on
$B_L=\pi^{-1}(L)$, there is a nonzero function~$h$ so that
$\omega = h\,\pi^*(dz)$ holds on~$B_L$.  Correspondingly, there will
be a function~$f$ on~$U$ so that 
$\bar a\,\bar\omega = \pi^*\left(f_{\bar z}\,d\bar z\right)$ and a
function~$g$ on~$U$ so that~$\bar s\,\omega^2 = \pi^*(g\,dz^2)$.  The
above structure equations then imply that the product~$e^{-3f}g$ is
holomorphic in~$z$.  Consequently, the quadratic form~$\bar s\,\omega^2$
is a nonvanishing multiple of a holomorphic quadratic form on~$L$ and
so either vanishes identically or else only vanishes at discrete points
of~$L$ and then only to finite order.  Note that~$|s|^2$
vanishes identically on a complex leaf if and only if that leaf is
totally geodesic in~$\bbP^2_R$.
\endremark

\remark{Remark 2}
It will be useful to understand the metric~$\omega\circ\bar\omega$ 
induced on the complex leaves, in particular, the Gauss curvature of 
this induced metric. Now, the equation for~$d\omega$ can be written
in the form
$$
d\omega = -i\rho\w\omega + 2s\,\bar\omega\w\eta
$$
where~$\rho=-3\tau + i(a\,\omega-\bar a\,\bar\omega)$. The equation
$$
d\rho \equiv -{\ts{\frac i2}}(4R-8|s|^2)\,\omega\w\bar\omega \bmod\eta
$$
then shows that the function~$K=4(R-2|s|^2)$ restricts to each complex
leaf to be its Gauss curvature. 
\endremark

\subsubhead 2.2.1. First case \endsubsubhead
Using the structure equations to expand the 
identity~$d(da)=0$ and then reducing the result 
modulo~$\omega$ yields
$$
s\,x - \bar s\,\bar x = 0.
$$

There are now two cases to consider.  First, suppose that~$s$ vanishes 
identically. Then so do~$x$ and $y$, and the remaining structure equation 
for~$a$ is
$$
da = -3ia\,\tau - 3a^2 \,\omega-({\ts\frac12}R{-}|a|^2)\,\bar\omega\,.
$$
Differentiating this equation just yields an identity.  Thus, the system
$$
\aligned
d\tau &= iR\,\omega\w\bar\omega\\
d\eta &= (a\,\omega + \bar a\,\bar\omega)\w\eta\\
d\omega &= (3i\tau - \bar a\,\bar\omega)\w\omega \\
da &= -3ia\,\tau - 3a^2 \,\omega-({\ts\frac12}R{-}|a|^2)\,\bar\omega\\
\endaligned
\tag1
$$
is differentially closed%
\footnote{I.e., the exterior derivatives of these equations are 
identities.  Of course, it then follows from Cartan's generalization 
of Lie's Third Fundamental Theorem that there are solutions to these 
equations, but the explicit computations in the next section will make
recourse to Cartan's theorem unnecessary.  This same comment applies
to the other two cases that will turn up in the next subsubsection.
} 
and describes the class of solutions~$\Sigma$ for 
which the complex leaves are totally geodesic.  This class will
be analyzed in the next section, after all of the integrability conditions 
have been found for the remaining cases.  

\subsubhead 2.2.2. Second and third cases \endsubsubhead
Suppose now that~$s$ does not vanish identically.  
Since~$\Sigma$ is real analytic and connected and since~$|s|^2$ is 
well-defined on~$\Sigma$,
there is a dense open set~$\Sigma^*\subset\Sigma$ 
on which~$|s|^2>0$.  On the bundle~$B_1^*=\pi^{-1}(\Sigma^*)\cap B_1$, which 
is a dense open subset of~$B_1$, write~$x = \bar s p$, where $p$ is real.  
The structure equations are now
$$
\align
da &= -3ia\,\tau-6\bar a\bar s\,\eta + (\bar s\,p-3a^2)\,\omega
            -({\ts\frac12}R{-}|a|^2{+}2|s|^2)\,\bar\omega\\
ds &= \phantom{-}6is\,\tau+ sp\,\eta + 3sa\,\omega + sy\,\bar\omega\\
\endalign
$$
(where, to simplify equations to follow, I have replaced the former~$y$ 
by~$sy$, which is permissible since~$s$ is nonzero). 

Now,~$a$ cannot vanish identically.  If it were to do so, 
then the above equations would imply~$p=0$ and $R = -4|s|^2 <0$ 
(since~$s$ is nonzero).  The equation for~$ds$ would then simplify 
to~$ds = 6is\,\tau + sy\,\bar\omega$.   Differentiating 
the relation~$R+4|s|^2=0$ then shows that~$y=0$, in turn implying 
that~$ds = 6is\,\tau$, which then implies that~$d\tau=0$, 
contradicting the structure equation for~$\tau$ since~$R\not=0$.  
By the real analyticity and connectedness of~$\Sigma$, it follows
that~$|a|^2$ is nonzero on a dense open set~$\Sigma^{**}\subset\Sigma^*$ 
and I can restrict attention to the corresponding subbundle~$B_1^{**}$, 
which I will do from now on.  Thus,~$a$ is nonzero on~$B_1^{**}$.

Now, the structure equations plus the reality of~$p$ yield
$$
0 = \frac{d\bigl(da\bigr)\w\bar\omega}{\bar s}
    + \frac{d\bigl(d\bar a\bigr)\w\omega}s
 = 6(\bar a \, \bar y - a y)\,\eta\w\omega\w\bar\omega.
$$
Thus~$ay$ is real, implying that there exists a function~$q=\bar q$ 
for which~$y = \bar a\,(q+3)$.
(Writing~$q{+}3$ instead of~$q$ here simplifies the following formulae.)  
Expanding the identity $d(da)=0$ and using the reality of~$p$ 
implies that~$p$ satisfies the equation
$$
\align
dp &= (2R{-}64|a|^2{+}8|s|^2{-}6|a|^2q{-}p^2)\eta\\
&\qquad\qquad-(ap + 24\bar a \bar s + 2\bar a \bar s q)\omega
    -(\bar a p + 24 as  + 2as q)\bar\omega\,.\\
\endalign
$$

By this structure equation and the reality of~$q$,
$$
0 = \frac{d\bigl(ds\bigr)\w a\,\omega}s
    + \frac{d\bigl(d\bar s\bigr)\w \bar a\,\bar\omega}{\bar s}
 = 4\,q\,(a^2s - {\bar a}^2\,\bar s)\,\eta\w\omega\w\bar\omega.
$$
Thus, either~$q$ or the imaginary part of~$a^2s$ vanishes
identically.  These two cases will be considered separately.

First, suppose that~$a^2s$ is real and introduce a real-valued
function~$t=\bar t$ so that~$s = {\bar a}^2\,t$. Using the reality of~$t$ 
and expanding the identities $0 = d(da)=d(ds)=d(dp)$ yields
$$
q = R + 4|a|^2+2|a|^2pt+|a|^4t^2
$$
plus a formula for~$dt$.  The result is structure equations of the form 
$$
\aligned
d\tau &= iR\,\omega\w\bar\omega\\
d\eta &= (a\,\omega + \bar a\,\bar\omega)\w\eta\\
d\omega &= (3i\tau - \bar a\,\bar\omega)\w\omega 
              + 2\,{\bar a}^2\,t\,\bar\omega\w\eta\\
da &= -3ia\,\tau-6|a|^2at\,\eta + a^2(t\,p-3)\,\omega
            -({\ts\frac12}R{-}|a|^2{+}2|a|^4t^2)\,\bar\omega\\
dp &= -(4R + 16|a|^2 + 16|a|^4t^2 + 12|a|^2pt + p^2)\eta\\
   &\qquad-\bigl(p + 8|a|^2t + 2t(R+4|a|^4t^2+2|a|^2pt)\bigr)
            (a\,\omega+\bar a\,\bar\omega)\\
dt &= t(p+12|a|^2t)\eta 
         + t(1+4|a|^2t^2+R/|a|^2)(a\,\omega+\bar a\,\bar\omega).\\
\endaligned
\tag2
$$
Differentiating these equations yields only identities, so this represents 
a set of solutions.  These will be analyzed below.  
This system is compatible with the relation~$t=0$,
in which case the structure equations specialize to~$(1)$, the first solution 
found.  Thus, the solutions~$(1)$ can be regarded as special cases of~$(2)$.

On the other hand, if~$q\equiv 0$, then the structure equations 
yield~$d\bigl(d(s)\bigr)=6sR\,\omega\w\bar\omega$, so this case can only 
occur when~$R=0$.  Assuming this, the structure equations found so far are
$$
\aligned
d\tau &= 0\\
d\eta &= (a\,\omega + \bar a\,\bar\omega)\w\eta\\
d\omega &= (3i\tau - \bar a\,\bar\omega)\w\omega + 2\,s\,\bar\omega\w\eta\\
da &= -3ia\,\tau-6\bar a\bar s\,\eta + (\bar s\,p-3a^2)\,\omega
            +\bigl(|a|^2{-}2|s|^2\bigr)\,\bar\omega\\
ds &= s\,(6i\tau + p\,\eta + 3a\,\omega + 3\bar a\,\bar\omega)\\
dp &= \bigl(8|s|^2{-}64|a|^2{-}p^2\bigr)\,\eta
        -(ap+24\bar a\,\bar s)\omega-(\bar a\,p+24 as)\bar\omega\,.\\
\endaligned
\tag3
$$
Differentiating these equations yield only identities, so this represents a
class of solutions that exist only in the case~$R=0$.  These will be analyzed
below.  Since~$a^2s$ is not, in general, real for these solutions,
they are not special cases of~$(2)$, although when~$s=0$, these solutions
do specialize to the~$t=0$ solutions of~$(2)$ in the case~$R=0$.  These
special solutions are the only overlap between the two.

\head 3. Existence of Solutions \endhead

In this section, I will prove general existence results that
assure that there are solutions to the equations~$(1)$, $(2)$, and $(3)$.
In each case, this will be followed by an analysis of the equations that
allows a complete description of the corresponding solutions.

\subhead 3.1. Solutions of type~1 \endsubhead

\subsubhead 3.1.1. Existence via the Frobenius theorem \endsubsubhead
Let~$M^{10} = G_R\times\bbC$ and let~$\gbold:M\to G_R$ and $\abold:M\to\bbC$ 
be the projections onto the factors.  I will regard forms on~$G_R$ or~$\bbC$
as forms on~$M$ via the pullbacks under these two maps and will not 
notate the pullback explicitly. Let~${\Cal I}_1$ be the 
exterior ideal on~$M$ generated by the linearly independent
real-valued $1$-forms~$\theta_1,\ldots,\theta_6$ where
$$
\aligned
            \theta_1 &= i(\bar\eta-\eta)\\
            \theta_2 &= \phi+i\,\abold\,\omega-i\,\bar\abold\bar\omega\\
\theta_3+i\,\theta_4 &= \sigma+2\bar\abold\,\eta \\
\theta_5+i\,\theta_6 &= d\abold+3i\,\abold\,\tau + 3\abold^2\,\omega
                          +({\ts{\frac12}}R-|\abold|^2)\,\bar\omega.\\
\endaligned
$$
The structure equations~$d\gamma =- \gamma\w\gamma$ imply that~${\Cal I}_1$ 
is differentially closed.  Thus, the Frobenius theorem implies that~$M$ 
is foliated by $4$-dimensional integral manifolds of~${\Cal I}_1$. 
Each leaf~$L\subset M$ is the image of a bundle
$B_1\subset G_R$ of a minimal Levi-flat hypersurface~$\Sigma$ satisfying 
equations~$(1)$ under the embedding~$\text{id}\times a:B_1\to G_R\times\bbC$. 
This gives an abstract description of the solutions of type~$(1)$.

Since~$G_R$ acts by left translation on~$G_R\times\bbC$ preserving the
ideal~${\Cal I}_1$, this left action permutes the integral manifolds,
and two integral manifolds are equivalent under this action if and only
if they correspond to congruent hypersurfaces in~$\bbP^2_R$. In particular,
two leaves~$L_1$ and $L_2$ represent equivalent solutions if and only if
they satisfy~$\abold(L_1)=\abold(L_2)$.  Note that this happens if and only
if the two images~$\abold(L_1)$ and $\abold(L_2)$ have nonempty 
intersection.

\subsubhead 3.1.2. Explicit description of the solutions \endsubsubhead
On any connected solution to~$(1)$, the structure equations imply
$$
\align
4\,da\w d\bar a 
&= \bigl(R+4|a|^2\bigr)
       \bigl( (R-8|a|^2)\,\omega\w\bar\omega
              +6i\,\tau\w(a\,\omega+\bar a\,\bar\omega\bigr) \bigr)\\
d(R+4|a|^2) &= -2\bigl(R+4|a|^2\bigr)(a\,\omega+\bar a\,\bar\omega).
\endalign
$$
It follows that for any leaf~$L$ of~${\Cal I}_1$, 
either the function~$R+4|\abold|^2$ vanishes identically or
else~$\abold:L\to\bbC$ is an immersion. 

Now, when~$R>0$, the only possibility is that~$\abold:L\to\bbC$ is
an immersion everywhere.  Moreover, using the left action of~$G_R$ plus
the existence of a leaf through any point of~$G_R\times\bbC$, it follows
that $\abold:L\to\bbC$ is a surjective submersion for every leaf.  In
particular, all of the leaves are equivalent under the action of~$G_R$. 
Since~$|s|^2$ vanishes identically on~$L$, it follows that the complex
leaves of~$\Sigma$ are totally geodesic in~$\bbP^2_R$, which is,
up to a constant scale factor, isometric to~$\bbC\bbP^2$ endowed with
the Fubini-Study metric.  Thus, $\Sigma$ must be a 1-parameter
family of complex lines in~$\bbC\bbP^2$.  In fact,~$\Sigma$ must be 
congruent to the smooth locus~$C^*_1$ of the `cone'
$$
C_1 
= \left\{\ \left[\matrix z\\ w\\ e^{ir}w\\ \endmatrix\right]\in\bbC\bbP^2
         \ \vrule\ r\in\bbR,\, [z,w]\in\bbC\bbP^1\ \right\}.
$$
It is evident that~$C^*_1$ is both Levi-flat and minimal.  Note that
$C_1$ has only one singular point (the intersection of the
complex lines that foliate it) and is otherwise smooth.

When~$R=0$, so that~$\bbP^2_0$ is isometric to~$\bbC^2$ with the
standard flat metric, there are two possibilities.  
The first possibility is that $|a|^2$ vanishes identically, in which case 
the corresponding~$\Sigma$ is congruent to a real hyperplane:
$$
H_0 
= \left\{\ \left[\matrix 1\\ z\\ r\\ \endmatrix\right]\in\bbP^2_0
         \ \vrule\ r\in\bbR,\ z\in\bbC\ \right\}.
$$
The second possibility is that~$|a|^2$ never vanishes.  By the same 
sort of argument made for the case of positive holomorphic sectional 
curvature, one sees that all of these cases are equivalent to the
smooth part of the cone
$$
C_0 
= \left\{\ \left[\matrix 1\\ z\\ e^{ir}\,z\\ \endmatrix\right]\in\bbP^2_0
         \ \vrule\ r\in\bbR,\,z\in\bbC\ \right\}.
$$

When~$R<0$, there is no loss of generality in setting~$R=-1$,
so I will do so for this discussion. Then
$$
\bbP^2_{-1} 
= \left\{\ \left[\matrix 1\\ z^1\\ z^2\\ \endmatrix\right]\in\bbC\bbP^2
         \ \vrule\ \bigl|z^1\bigr|^2+\bigl|z^2\bigr|^2 <1\ \right\}
$$
is the hyperbolic complex $2$-ball and there are three possibilities, 
depending on the sign of~$R+4|a|^2=4|a|^2-1$.

The solutions with~$4|a|^2-1>0$ are all congruent to the smooth part of
the hyperbolic version of the cone:
$$
C_{-1} 
= \left\{\ \left[\matrix 1\\ z\\ e^{ir}\,z\\ \endmatrix\right]\in\bbP^2_{-1}
         \ \vrule\ r\in\bbR,\ \sqrt2\,|z|<1\ \right\}.
$$
This cone has one singular point.  The leaves of~$dr=0$ are the complex
leaves, each one biholomorphic to a punctured disk.

All solutions with~$4|a|^2-1=0$ are congruent to the `horosphere' solution
$$
S_{-1} 
= \left\{\ \left[\matrix 1\\ z\\ ir(1-z)\\ \endmatrix\right]\in\bbP^2_{-1}
      \ \vrule\ r\in\bbR,\ |z|^2+r^2|1{-}z|^2 < 1 \ \right\}.
$$
The complex leaves in~$S_{-1}$ are the leaves of~$dr=0$ in 
the chosen parametrization.  All of these complex leaves
intersect at one point on the boundary of the ball.  This solution can be
interpreted as a limit of the cone~$C_{-1}$ as one moves the singular point 
of the cone out to the boundary of~$\bbP^2_R$ in~$\bbP^2$.

All solutions with~$4|a|^2-1<0$ are congruent to
the hyperbolic version of the hyperplane solution, namely
$$
H_{-1} 
= \left\{\ \left[\matrix 1\\ z\\ r\\ \endmatrix\right]\in\bbP^2_0
         \ \vrule\ r\in\bbR,\ r^2+|z|^2<1\ \right\}.
$$

This completes the list of solutions of the system~$(1)$.

\subhead 3.2. Solutions of type~2 \endsubhead
Consider the solutions of the system~$(2)$.  To avoid repetition, I am
going to consider only solutions for which~$t$ is non-zero, since the
solutions with $t$ vanishing identically have already been accounted
for as solutions of type~$(1)$.

\subsubhead 3.2.1. Existence via the Frobenius theorem \endsubsubhead
Let $M^{12}=G_R\times \bbC^*\times\bbR\times\bbR$, and let $\gbold:M\to G_R$,
$\abold:M\to\bbC^*$, $\pbold:M\to\bbR$, and $\tbold:M\to\bbR$ be the
projections onto the first through fourth factors, respectively.  
Let~${\Cal I}_2$ be the exterior ideal on~$M$ generated by the linearly 
independent real-valued $1$-forms~$\theta_1,\ldots,\theta_8$ where
$$
\aligned
            \theta_1 &= i(\bar\eta-\eta)\\
            \theta_2 &= \phi+i\,\abold\,\omega-i\,\bar\abold\bar\omega\\
\theta_3+i\,\theta_4 &= \sigma+2\bar\abold\,\eta 
                           +2\bar\abold^2\tbold\,\bar\omega\\
\theta_5+i\,\theta_6 &= d\abold+3i\,\abold\,\tau 
                      +6|\abold|^2\abold\tbold\,\eta 
                    - \abold^2(\tbold\,\pbold-3)\,\omega
        +({\ts\frac12}R{-}|\abold|^2{+}2|\abold|^4\tbold^2)\,\bar\omega\\
\theta_7 &= 
      d\pbold + 4(R + 4|\abold|^2 + 4|\abold|^4\tbold^2 
+ 3|\abold|^2\pbold\tbold + \pbold^2)\eta\\
         &\qquad + \bigl(\pbold 
+ 2\tbold(R+4|\abold|^2+4|\abold|^4\tbold^2+2|\abold|^2\pbold\tbold)\bigr)
            (\abold\omega+\bar \abold\,\bar\omega)\\
\theta_8 &= d\tbold 
- \tbold\bigl((\pbold+12|\abold|^2\tbold)\eta 
              + (1+4|\abold|^2\tbold^2+R/|\abold|^2)
                (\abold\omega+\bar \abold\,\bar\omega)\bigr).\\
\endaligned
$$
(The reason for the restriction~$\abold\not=0$ is the division by~$|\abold|^2$
in the last formula.)   The structure equations show that 
the ideal~${\Cal I}_2$ is closed under exterior differentiation, so~$M$
is foliated by $4$-dimensional integral manifolds of~${\Cal I}_2$.  

By construction, each leaf~$L\subset M$ is the image of the bundle
$B^{**}_1\subset G_R$ over the nondegenerate
part~$\Sigma^{**}$ of a minimal Levi-flat hypersurface~$\Sigma$ satisfying 
equations~$(2)$ under the embedding~
$$
\text{id}\times a\times p\times t:B_1
\longrightarrow G_R\times\bbC^*\times\bbR\times\bbR = M.
$$

Since~$G_R$ acts by left translation on~$G_R\times\bbC^*\times\bbR\times\bbR$ 
preserving the ideal~${\Cal I}_2$, this left action permutes its integral 
manifolds, and two integral manifolds are equivalent under this action if and 
only if they correspond to congruent hypersurfaces in~$\bbP^2_R$. 
In particular, two leaves~$L_1$ and $L_2$ represent equivalent solutions 
if and only if they satisfy~$(\abold,\pbold,\tbold)(L_1)
=(\abold,\pbold,\tbold)(L_2)$.  

In fact, in order for two leaves~$L_1$
and~$L_2$ to be equivalent under~$G_R$, it suffices that the two image sets
$(\abold,\pbold,\tbold)(L_1)$ and $(\abold,\pbold,\tbold)(L_2)$ in
$\bbC^*\times\bbR\times\bbR$ have a nonempty intersection.  To see why this
is so, note that if~$L_i$ contains~$(g_i,a,p,t)$, then the 
submanifold~$L$ described by
$$
L = \left\{\,(g_2{g_1}^{-1}g,b,q,u)\,\vrule\, (g,b,q,u)\in L_1\right\}
$$
contains~$(g_2,a,p,t)\in L_2$, is evidently a maximal integral manifold 
of~${\Cal I}_2$, and so must equal~$L_2$.  In particular, in order to
classify the solutions up to rigid motion, it would suffice to determine
the partition of~$\bbC^*\times\bbR\times\bbR$ into the images
$(\abold,\pbold,\tbold)(L)$ as~$L$ ranges over the leaves of~${\Cal I}_2$.
Moreover, this argument shows that the fibers of the 
map~$(\abold,\pbold,\tbold):L\to\bbC^*\times\bbR\times\bbR$ are the 
orbits of the action on~$L$ of the ambient symmetry group 
of the corresponding solution~$\Sigma^{**}$. 

The structure equations imply that the function~$\tbold$ cannot vanish
anywhere on a leaf~$L$ unless it vanishes identically on~$L$.  As mentioned
at the begining of this subsection, the leaves on which~$\tbold$ vanishes
identically are of type~$(1)$ and so can be set aside in this discussion.
For the rest of this subsection, the assumption that~$\tbold$ is nonvanishing
on~$L$ will be in force.

\subsubhead 3.2.2. Symmetries of the solutions \endsubsubhead
One might expect the images~$(\abold,\pbold,\tbold)(L)$
to have dimension~$4$, at least at `generic' points, since each leaf~$L$
has dimension~$4$.  However, equations~$(2)$ imply 
that~$da\w d\bar a\w dp\w dt$ vanishes identically.  Consequently, the rank 
of~$(\abold,\pbold,\tbold):L\to\bbC^*\times\bbR\times\bbR$ is strictly
less than~$4$ at all points, implying that the fibers of this map 
(and hence the symmetry group of~$L$) must have positive dimension.

It is not hard to make these fibers explicit. 
By the structure equations~$(2)$, the (real) nowhere vanishing vector 
field~$Y$ on~$L$ that satisfies
$$
\align
\tau(Y) &= R + 2|a|^2(pt-4) + 4|a|^4t^2\\
\eta(Y) &=0\\
\omega(Y) &= 6i\,\bar a\\
\endalign
$$
also satisfies~$da(Y)=dp(Y)=dt(Y)=0$. The structure equations also
show that, for the generic value
$(a_0,p_0,t_0)\in\bbC^*{\times}\bbR{\times}\bbR$, the leaf~$L$ whose
$(\abold,\pbold,\tbold)$@-image contains~$(a_0,p_0,t_0)$ has the 
property that~$(d\abold,d\pbold,d\tbold)$ has rank~$3$ along the
preimage of~$(a_0,p_0,t_0)$.  In particular,~$Y$ spans the tangent
to the fiber at such points.  

Given this, it would not be surprising to find that~$Y$ can be
scaled so as to become a symmetry vector field.  In fact, one finds
that the flow of~$X = e^f\,Y$ preserves the 
coframing~$(\tau,\rho,\omega)$ on~$L$ if and only if~$f$ satisfies the
equation
$$
df = 8|a|^2\,t\,\eta - (pt-2)(a\,\omega+\bar a\,\bar\omega)
$$
Now, by the structure equations, the right hand side of this equation 
is a closed $1$-form on~$L$.  This shows that, at least locally
(or, more precisely, on some covering space of~$L$), a scaling factor~$e^f$
exists making~$X=e^f\,Y$ a symmetry vector field.  Moreover, this~$f$
is unique up to the addition of a constant.

This implies that any solution hypersurface~$\Sigma\subset\bbP^2_R$ whose 
structure equations are of the form~$(2)$ must actually be invariant under 
a one-parameter group of isometries of~$\bbP^2_R$, i.e., a one-parameter
subgroup of~$G_R$.  Moreover, because~$\eta(Y)=0$, this one-parameter 
subgroup can be chosen (if it is not actually unique) so that it preserves 
each complex leaf in~$\Sigma$.   

\subsubhead 3.2.3. Explicit solutions invariant under a given 1-parameter
subgroup \endsubsubhead
A one-parameter subgroup of isometries of~$\bbP^2_R$ is of the 
form~$\{e^{tz}\,\vrule\,t\in\bbR\}$ for some~$z\not=0$ in~$\eug_R$.  
There is a unique holomorphic vector field~$Z$
on~$\bbP^2_R$ whose real part is the infinitesimal generator of the action 
of the subgroup~$\{e^{tz}\,\vrule\,t\in\bbR\}$.  If~$U_z\subset\bbP^2_R$ 
denotes the open set that is the complement of the fixed locus of the 
flow~$e^{tz}$, then~$U_z$ is foliated by complex curves that are the 
`integral curves' of the holomorphic flow generated by~$Z$.  
By the above discussion, the nondegenerate part~$\Sigma^{**}$ of any 
solution~$\Sigma$ of type~$(2)$ will be swept out by a (real) one-parameter 
family of integral curves of~$Z$ for some isometric flow~$e^{tz}$.  Since,
by construction, the complex leaves of a solution of type~$(2)$ are not
totally geodesic, this shows that~$z\in\eug_R$ must be chosen so that
the $Z$-integral curves in~$U_z$ are not totally geodesic.  I will
refer to a~$z\in\eug_R$ with this property as~{\it nondegenerate\/}. 

Conversely, starting with any one-parameter subgroup~$e^{tz}$ of isometries
of~$\bbP^2_R$ and considering the corresponding holomorphic 
foliation of~$U_z\subset\bbP^2_R$ by complex curves, one can construct 
$e^{tz}$@-invariant Levi-flat hypersurfaces in~$U_z$ by taking the union of 
any (real) one-parameter family of complex leaves of this foliation.  
It now suffices to show that one can choose this one-parameter family in 
such a way that the resulting hypersurface will be minimal.  
I am going to show that this can always be done,
essentially in only one way up to isometry, and that, when~$z$ is 
nondegenerate in the sense of the previous paragraph, this always yields a 
solution~$\Sigma$ of type~$(2)$.  Thus, the solutions of type~$(2)$ 
correspond to the conjugacy classes of nondegenerate one-parameter subgroups 
of isometries of~$\bbP^2_R$.

First, consider the case where~$R>0$.  Without essential loss of generality,
I can assume that~$R=1$, so that~$G_R=G_1=\SU(3)$.  Every one-parameter
subgroup of~$\SU(3)$ is semi-simple and hence conjugate to a diagonal
subgroup generated by a nonzero element
$$
z = \pmatrix i\lambda_0&0&0\\ 0&i\lambda_1&0\\ 0&0&i\lambda_2\\ \endpmatrix
\qquad\text{where}\qquad \lambda_0+\lambda_1+\lambda_2=0.
$$
The corresponding vector field on~$\bbC^3$ (which is also well-defined on
$\bbP^2_R\simeq\bbP^2$) can be written in terms of unitary holomorphic
coordinates~$z = (z^a)$ as the real part of the holomorphic vector field
$$
Z =  i\lambda_0\,z^0\,\frac{\partial\hfil}{\partial z^0}
    +i\lambda_1\,z^1\,\frac{\partial\hfil}{\partial z^1}
    +i\lambda_2\,z^2\,\frac{\partial\hfil}{\partial z^2}.   
$$
The holomorphic integral curve of~$Z$ through~$c = [c^a]\in\bbP^2_R$ 
is of the form
$$
\left\{\ \left[\matrix c^0\,e^{i\lambda_0w}\\
                c^1\,e^{i\lambda_1w}\\
                c^2\,e^{i\lambda_2w}\endmatrix\right]
       \ \vrule\ w\in\bbC\/\right\}.
$$
This will be a point or a line for all such~$c$ if and only if
two of the $\lambda_i$ are equal.  In such a case, the integral curves 
of~$Z$ are open subsets of lines through a fixed point in~$\bbP^2$.  
Thus any minimal Levi-flat hypersurface whose 
complex leaves are integral curves of~$Z$ will be of type~$(1)$.  
Set this case aside and, from now on, assume that~$z$ is nondegenerate, 
i.e., that the $\lambda_i$ are mutually distinct. 

Since radial dilation has no effect on the projective space, the flow
the vector field~$Z$ induces on~$\bbP^2$ is the same as that of the
vector field
$$
Z' = i(\lambda_1-\lambda_0)\,z^1\,\frac{\partial\hfil}{\partial z^1}
    +i(\lambda_2-\lambda_0)\,z^2\,\frac{\partial\hfil}{\partial z^2}  
$$
and this, in turn, will have the same holomorphic integral curves
in~$\bbP^2$ as
$$
Z'' = z^1\,\frac{\partial\hfil}{\partial z^1}
    +\lambda\,z^2\,\frac{\partial\hfil}{\partial z^2} 
\qquad\text{where}\quad 
\lambda = \frac{(\lambda_2-\lambda_0)}{(\lambda_1-\lambda_0)}\not=0,1.
$$
The nonlinear integral curves of this vector field are of the form
$$
\left\{\ \left[\matrix 1\\
                c\,e^{w}\\
                e^{\lambda w}\endmatrix\right]
       \ \vrule\ w\in\bbC\/\right\}.
$$
where~$c$ is any nonzero complex constant.  Thus, a Levi-flat hypersurface
whose complex leaves are integral curves of this vector field can be
locally parametrized in the form
$$
\Sigma = \left\{\ \left[\matrix 1\\
                e^{w+x(r)+iy(r)}\\
                e^{\lambda w}\endmatrix\right]
       \ \vrule\ w\in\bbC, r\in I\/\right\}.
$$
where~$x+i\,y:I\to\bbC$ is some smooth immersion of an interval~$I\subset\bbR$.
Brute force calculation then yields that such a hypersurface is minimal if 
and only if~$y$ is a constant function.  Thus, up to a holomorphic isometry, 
such a minimal Levi-flat hypersurface is an open subset of the hypersurface
$$
\Sigma_\lambda = \left\{\ \left[\matrix 1\\
                e^{w+r}\\
                e^{\lambda w}\endmatrix\right]
       \ \vrule\ w\in\bbC, r\in \bbR\/\right\}.
$$
Note that~$\Sigma_\lambda$ is congruent to~$\Sigma_{1/\lambda}$ but that, 
otherwise, the~$\Sigma_\lambda$ are mutually noncongruent.
When~$\lambda$ is irrational, this hypersurface is dense in~$\bbP^2$, 
but when~$\lambda = p/q$ where~$p$ ($\not=0,q$) and $q>0$ are integers without
common factors, this hypersurface is dense in an algebraically defined
hypersurface that is singular at the point~$z^1=z^2=0$ but can also
be singular along the entire lines~$z^1=0$ and $z^2=0$, depending on
the values of~$p$ and~$q$.  A typical such hypersurface is defined
by an equation of the 
form~$\text{Im}\bigl(({\bar z}^0)^{p+q}(z^1)^{p}(z^2)^{q}\bigr)=0$.

Next, consider the case where~$R=0$, i.e., when~$\bbP^2_R=\bbP^2_0$ is 
isometric to~$\bbC^2$ with its standard flat metric.  Let~$(z^1,z^2)$ 
be unitary holomorphic linear coordinates on~$\bbC^2$.  A nonzero vector 
field whose flow is a holomorphic isometry on~$\bbC^2$ is then
conjugate via an action of~$G_0$ to a constant multiple of the real part 
of either
$$
Z =  i\,z^1\,\frac{\partial\hfil}{\partial z^1}
    +i\lambda\,z^2\,\frac{\partial\hfil}{\partial z^2} 
\qquad\text{or}\qquad
Z =  i\,\frac{\partial\hfil}{\partial z^1}
    +i\lambda\,z^2\,\frac{\partial\hfil}{\partial z^2} 
$$
for some real number~$\lambda$.  In the first case, the holomorphic integral
curves of~$Z$ will be lines in~$\bbC^2$ if and only if~$\lambda=0$~or~$1$ 
while, in the second case, the holomorphic integral curves of~$Z$ will be 
lines in~$\bbC^2$ if and only if~$\lambda=0$.  These are the degenerate
values that will be set aside, as these degenerate cases lead to the
hyperplane or Clifford cone solutions that have already been discussed
in the previous subsection.

Consider the first type of vector field with~$\lambda\not=0$ or~$1$.  Any
holomorphic integral curve of~$Z$ that is not contained in a line in~$\bbC^2$
is of the form
$$
\left\{\pmatrix c\,e^{w}\\
                \,e^{\lambda w}\endpmatrix
       \ \vrule\ w\in\bbC\/\right\}.
$$
where~$c\in\bbC$ is a nonzero constant.   A smooth Levi-flat 
hypersurface~$\Sigma^3\subset\bbC^2$ whose complex leaves consist of 
such integral curves can be locally parametrized in the form
$$
\Sigma^3 = \left\{\pmatrix e^{w+x(r)+iy(r)}\\
                \,e^{\lambda w}\endpmatrix
       \ \vrule\ w\in\bbC, r\in I\/\right\}.
$$
where~$x+i\,y:I\to\bbC$ is some smooth immersion of an interval~$I\subset\bbR$.
Brute force calculation then yields that such a hypersurface is minimal if and
only if~$y$ is a constant function.  Consequently, it follows that, up to
a holomorphic isometry, the connected solutions of this kind are all 
equivalent to open subsets of the immersed hypersurface
$$
\Sigma_\lambda = \left\{\pmatrix e^{w+r}\\
                   \,e^{\lambda w}\endpmatrix
                   \ \vrule\ w\in\bbC, r\in \bbR\/\right\}.
$$
If $\lambda$ is irrational, then $\Sigma_\lambda$ is dense in~$\bbC^2$ and
the (implicitly described) immersion given above is an embedding.  On the
other hand, if~$\lambda = p/q$ where~$p\not=0$ and $q>0$ are distinct integers 
without common factors, then this immersion is not an embedding.  Moreover,
$\Sigma_{p/q}$ is dense in an algebraic real hypersurface, namely
$$
\align
\left(z^1\right)^p\left({\bar z}^2\right)^q 
- \left({\bar z}^1\right)^p\left(z^2\right)^q = 0 
&\quad \text{when $p>0$},\\
\left({\bar z}^1\right)^{-p}\left({\bar z}^2\right)^q 
- \left(z^1\right)^{-p}\left(z^2\right)^q = 0 
&\quad \text{when $p<0$}.\\
\endalign
$$
Note that these hypersurfaces are cones that are singular at the origin
and along the axes except when $p$ or~$q$ equals~$1$.

Consider the second type of vector field with~$\lambda\not=0$.  Any
holomorphic integral curve of~$Z$ that is not contained in a line in~$\bbC^2$
is of the form
$$
\left\{\pmatrix w+c\\
                \,e^{\lambda w}\endpmatrix
       \ \vrule\ w\in\bbC\/\right\}.
$$
where~$c\in\bbC$ is a nonzero constant.   A smooth Levi-flat 
hypersurface~$\Sigma^3\subset\bbC^2$ whose complex leaves consist of 
such integral curves can be locally parametrized in the form
$$
\Sigma^3 = \left\{\pmatrix w+x(r)+iy(r)\\
                \,e^{\lambda w}\endpmatrix
       \ \vrule\ w\in\bbC, r\in I\/\right\}.
$$
where~$x+i\,y:I\to\bbC$ is some smooth immersion of an interval~$I\subset\bbR$.
Brute force calculation then yields that such a hypersurface is minimal if and 
only if~$y$ is a constant function.  Consequently, up to a holomorphic 
isometry followed by a homo\-thety, the connected solutions of this kind are 
open subsets of the closed embedded hypersurface
$$
\Sigma = \left\{\pmatrix w\\
                   \,re^w\endpmatrix
                   \ \vrule\ w\in\bbC, r\in \bbR\/\right\}.
$$
(In this parametrization, the complex leaf given by~$r=0$ does not belong
to~$\Sigma^*$, as defined in~\S2.2.2.)  This hypersurface can be defined
implicitly by the equation 
$$
\text{Im}(z^2 e^{-z^1})=0
$$
and is evidently transcendental.

Finally, consider the case~$R<0$, where, without essential loss of 
generality, it suffices to consider only the case~$R=-1$.  The equivalence
classes of one-dimensional subspaces of~$\eusu(2,1)=\eug_{-1}$ under the 
adjoint action are more complicated in this case.  The elements~$z$ that 
have an eigenvector that is $\la,\ra_{-1}$-positive (and whose 
associated flow, therefore, has a fixed point in~$\bbP^2_{-1}$) can
be diagonalized in the form
$$
z = \pmatrix i\lambda_0&0&0\\ 0&i\lambda_1&0\\ 0&0&i\lambda_2\\ \endpmatrix
\qquad\text{where}\qquad \lambda_0+\lambda_1+\lambda_2=0.
$$
If $z$ has no $\la,\ra_{-1}$-positive eigenvector, then it must have a
null eigenvector.  In this case, the most generic possibility is for~$z$
to have three distinct eigenvalues, in which case two of the eigenvalues
cannot be purely imaginary and their corresponding eigenvectors must be
$\la,\ra_{-1}$-null.  Consequently, one can normalize these eigenvectors
and show that, up to a (real) multiple,~$z$ is conjugate to an element of 
the form
$$
z = \pmatrix i\lambda&1&0\\ 1&i\lambda&0\\ 0&0&-2i\lambda\\ \endpmatrix
\qquad\text{where}\quad~\lambda\in\bbR.
$$
If~$z$ has a double eigenvalue with a unique 
$\la,\ra_{-1}$-null corresponding eigenvector 
and a $\la,\ra_{-1}$-negative eigenvector, then,
up to a (real) multiple, $z$ is conjugate to an element of the form
$$
z = \pmatrix i(\lambda{+}1)&-i&0\\ 
              i&i(\lambda{-}1)&0\\ 0&0&-2i\lambda\\ \endpmatrix
\qquad\text{where}\quad~0\not=\lambda\in\bbR.
$$
If~$z$ has a triple eigenvalue, i.e., is nilpotent, then either~$z^2\not=0$,
in which case it is conjugate to an element of the form
$$
z = \pmatrix 0&0&1\\ 0&0&1\\ 1&-1&0\\ \endpmatrix,
$$
or else~$z^2=0$ (the most degenerate case), in which case it is conjugate
to an element of the form
$$
z = \pmatrix i&-i&0\\ 
             i&-i&0\\ 0&0&0\\ \endpmatrix.
$$

Among these five cases, the holomorphic flow on~$\bbP^2_{-1}$ corresponding
to~$e^{tz}$ will have all integral curves be totally geodesic in two cases.  
In the case where~$z$ is diagonalizable, this happens 
when $\{\lambda_0,\lambda_1,\lambda_2\}$ are not distinct.  Among the 
nondiagonalizable cases, this happens only for the last case, i.e., when
$z^2 = 0$.  These cases will be set aside, as they have already been
treated in the discussion of type~$(1)$ solutions.

Now, in the diagonalizable case, the analysis proceeds exactly along
the lines of the elliptic case and there is no need to give details.  The 
end result is that a connected Levi-flat minimal hypersurface whose complex 
leaves are invariant under a nondegenerate diagonalizable flow is congruent 
to an open subset of the hypersurface
$$
\Sigma_\lambda = \left\{\ \left[\matrix 1\\
                e^{w+r}\\
                e^{\lambda w}\endmatrix\right]
       \ \vrule\ w\in\bbC, r\in \bbR, 
       \left| e^{w+r}\right|^2 + \left| e^{\lambda w}\right|^2 < 1 \/\right\}
\subset \bbP^2_{-1}.
$$
where~$\lambda$ is a real constant not equal to $0$ or $1$.  This hypersurface
has an algebraic defining equation if and only if~$\lambda$ is rational.

The next case, where~$z$ has two distinct $\la,\ra_{-1}$-null eigenvectors, 
can be analyzed in a similar manner and one finds that a connected Levi-flat
minimal hypersurface whose complex leaves are invariant under the 
associated holomorphic flow is congruent to an open subset of the 
hypersurface
$$
\Sigma'_\lambda = \left\{\ \left[\matrix 
                e^{(1+i\lambda)(w+r)}+e^{(-1+i\lambda)w}\\
                e^{(1+i\lambda)(w+r)}-e^{(-1+i\lambda)w}\\
                1\endmatrix\right]
       \ \vrule\ w\in\bbC, r\in \bbR 
       \/\right\}\subset \bbP^2_{-1}.
$$
where~$\lambda$ is a real constant.  When~$\lambda=0$, this is a real
curve in the pencil of conics that pass through two points on the 
boundary of~$\bbP^2_{-1}$ and have given tangents there.  When~$\lambda$
is nonzero, the curves~$r=r_0$ are not algebraic.  (Of course, $w$ and~$r$
must satisfy an inequality in order that the formula given in this description
represent at point in~$\bbP^2_{-1}$, but it is not useful to
make this inequality explicit for the purposes at hand.)

In the case where~$z$ has a double eigenvalue (and not a triple one),
a similar analysis shows that a connected Levi-flat minimal hypersurface
whose complex leaves are invariant under the associated holomorphic flow
is congruent to an open subset of the hypersurface
$$
\Sigma''_\mu = \left\{\ \left[\matrix 
                w+r+1\\
                w+r-1\\
                e^{\mu w}\endmatrix\right]
       \ \vrule\ w\in\bbC, r\in \bbR, \ 4r>e^{\mu(w+\bar{w})}-2(w{+}\bar{w})\ 
       \/\right\}\subset \bbP^2_{-1}.
$$
where~$\mu$ is a nonzero real constant.

In the final nondegenerate case, where the symmetry generator~$z\in\eug_{-1}$ 
satisfies $z^2\not=0$ but~$z^3 = 0$, the nonlinear integral curves of the 
associated holomorphic flow are conics (i.e., rational curves of degree~$2$) 
in~$\bbP^2$, all tangent at a point on the boundary of~$\bbP^2_{-1}\subset
\bbP^2$.  Brute force calculation shows that any Levi-flat minimal 
hypersurface whose complex leaves are invariant under such a flow is
congruent to the hypersurface
$$
\Sigma = \left\{\ \left[\matrix w^2+r+1\\
                w^2+r-1\\
                2w\endmatrix\right]
       \ \vrule\ w\in\bbC,\ r\in \bbR,\  
       4(\text{Im}\,w)^2 < r \/\right\}
\subset \bbP^2_{-1}.
$$
Details will be left to the reader.

\subhead 3.3. Solutions of type~$3$ \endsubhead
Finally, consider the solutions of the system~$(3)$.  To avoid repetition,
I will set aside the cases where the solution reduces to one of type~$(1)$.
This means that the solution has~$s\not=0$, which, by the structure 
equations~$(3)$, implies that~$s$ is nowhere vanishing.

\subsubhead 3.3.1. Existence via the Frobenius theorem \endsubsubhead
Let $M^{13}=G_0\times \bbC\times\bbC\times\bbR$, and let $\gbold:M\to G_0$,
$\abold:M\to\bbC$, $\sbold:M\to\bbC$, and $\pbold:M\to\bbR$ be the
projections onto the first through fourth factors, respectively.  
Let~${\Cal I}_3$ be the exterior ideal on~$M$ generated by the linearly 
independent real-valued $1$-forms~$\theta_1,\ldots,\theta_9$ where
$$
\aligned
            \theta_1 &= i(\bar\eta-\eta)\\
            \theta_2 &= \phi+i\,\abold\,\omega
                               -i\,\bar\abold\,\bar\omega\\
\theta_3+i\,\theta_4 &= \sigma+2\bar\abold\,\eta 
                           +2\sbold\,\bar\omega\\
\theta_5+i\,\theta_6 &= d\abold+3i\,\abold\,\tau 
                      +6\,\bar\abold\,\bar\sbold\,\eta 
                    - (\bar\sbold\,\pbold-3\abold^2)\,\omega
                    - (|\abold|^2{-}2|\sbold|^2)\,\bar\omega\\
\theta_7+i\,\theta_8 &= d\sbold 
           - \sbold\bigl( 6i\,\tau  + \pbold\,\eta 
                          + \abold\omega+ \bar\abold\,\bar\omega\bigr).\\
\theta_9 &= 
      d\pbold - (8|\sbold|^2-64|\abold|^2-\pbold^2)\eta
         + (\abold\,\pbold + 24\,\bar\abold\,\bar\sbold)\,\omega
                 + (\bar\abold\,\pbold + 24\,\abold\,\sbold)\,\bar\omega\,.\\
\endaligned
$$
By the structure equations, 
the ideal~${\Cal I}_3$ is closed under exterior differentiation, so~$M$
is foliated by $4$-dimensional integral manifolds of~${\Cal I}_3$.  

By construction, each leaf~$L\subset M$ is the image of the bundle
$B^{**}_1\subset G_0$ over the nondegenerate
part~$\Sigma^{**}$ of a minimal Levi-flat hypersurface~$\Sigma$ satisfying 
equations~$(3)$ under the embedding~
$$
\text{id}\times a\times s\times p:B^{**}_1
\longrightarrow G_0\times\bbC\times\bbC\times\bbR = M.
$$

Since~$G_0$ acts by left translation on~$G_0\times\bbC^*\times\bbR\times\bbR$ 
preserving the ideal~${\Cal I}_3$, this left action permutes its integral 
manifolds, and two integral manifolds are equivalent under this action if and 
only if they correspond to congruent hypersurfaces in~$\bbP^2_0\simeq\bbC^2$. 
In particular, two leaves~$L_1$ and $L_2$ represent equivalent solutions 
if and only if they satisfy~$(\abold,\sbold,\pbold)(L_1)
=(\abold,\sbold,\pbold)(L_2)$.  

In fact, in order for two leaves~$L_1$
and~$L_2$ to be equivalent under~$G_0$, it suffices that the two image sets
$(\abold,\sbold,\pbold)(L_1)$ and $(\abold,\sbold,\pbold)(L_2)$ in
$\bbC\times\bbC\times\bbR$ have a nonempty intersection.  To see why this
is so, note that if~$L_i$ contains~$(g_i,a,s,p)$, then the 
submanifold~$L$ described by
$$
L = \left\{\,(g_2{g_1}^{-1}g,b,q,u)\,\vrule\, (g,b,q,u)\in L_1\right\}
$$
contains~$(g_2,a,s,p)\in L_2$, is evidently a maximal integral manifold 
of~${\Cal I}_3$, and so must equal~$L_2$.  In particular, in order to
classify the solutions up to rigid motion, it would suffice to determine
the partition of~$\bbC\times\bbC\times\bbR$ into the images
$(\abold,\sbold,\pbold)(L)$ as~$L$ ranges over the leaves of~${\Cal I}_3$.

\subsubhead 3.3.2. First integrals and the symmetry of solutions
\endsubsubhead
Now, it would be reasonable to expect the images~$(\abold,\sbold,\pbold)(L)$
to have dimension~$4$, at least at `generic' points, since each leaf~$L$
has dimension~$4$.  In fact, by the argument in the previous paragraph, 
it is evident that the fibers of the map~$(\abold,\sbold,\pbold):L\to
\bbC\times\bbC\times\bbR$ are the orbits of the action on~$L$ of 
the ambient symmetry group of the corresponding solution~$\Sigma^{**}$. 

Consider the quantities%
\footnote{The significance of these quantities will become clear in
the analysis to be carried out below.}
$$
\aligned
\Abold &= {\ts{\frac19}}|\sbold|^{2/3}
                \bigl(48\,|\abold|^2+12\,|\sbold|^2+\pbold^2\bigr),\\
\Bbold &= {\ts{\frac1{27}}}|\sbold|\,
                \bigl(216\,\abold^2\sbold+216\,{\bar\abold}^2\bar\sbold
                      +72\,|\abold|^2\pbold-36\,|\sbold|^2\pbold+\pbold^3\bigr).\\ 
\endaligned
$$
The structure equations show that the 1-form~$d(\Abold^3{-}\Bbold^2)$ lies
in~${\Cal I}_3$, which implies that the image~$(\abold,\sbold,\pbold)(L)$ of
any ${\Cal I}_3$@-leaf~$L$ lies in a level set of~$\Fbold=\Abold^3{-}\Bbold^2$, 
a homogeneous polynomial of degree~$8$ in the variables~$\abold$, $\bar\abold$,
$\sbold$, $\bar\sbold$, and $\pbold$.

Calculation shows that~$\Fbold\ge0$, with 
equality exactly along the $3$@-dimensional 
cone~$C_0\subset\bbC\times\bbC\times\bbR$ defined by the equations
$$
0=\abold^2\sbold-\bar\abold^2\,\bar\sbold
 = 8\,|\abold|^2|\sbold|^2-4\,|\sbold|^4
     -2\,\abold^2\sbold\pbold-2\,\bar\abold^2\bar\sbold\pbold 
     +\pbold^2\,|\sbold|^2.
$$  
In particular, the ${\Cal I}_3$@-leaves that lie in~$G_0\times C_0$
represent either solutions of type~$(1)$ or of type~$(2)$,
and have already been analysed in the previous subsections.  
Moreover,~$0$ is the only critical value 
of~$\Fbold$ on~$\bbC\times\bbC\times\bbR\simeq\bbR^5$. 
The remaining level sets of~$\Fbold$ are smooth, connected hypersurfaces.  
In fact, because~$\Fbold$ is a homogeneous polynomial of degree~$8$, 
it follows that all of the positive level sets are diffeomorphic 
by homothety.  

A rather laborious calculation using the structure equations above shows that 
for any ${\Cal I}_3$@-leaf~$L$ on which~$\Fbold=c^2>0$, the rank of the 
map~$(\abold,\sbold,\pbold):L\to\bbC\times\bbC\times\bbR$ is~$4$, i.e.,
that~$(\abold,\sbold,\pbold):L\to\Fbold^{-1}(c^2)$ is a local diffeomorphism.
The existence theorem proved above via the Frobenius theorem coupled with
the $G_0$-invariance of~${\Cal I}_3$ shows that this map must actually be
a (surjective) covering map.  Thus, there is a 1-parameter family of 
noncongruent solutions of type~$(3)$, one for each positive level set 
of~$\Fbold$.  

\subsubhead 3.3.3. The effect of homothety
\endsubsubhead
While the members of this 1-parameter family are mutually incongruent by
isometries, it turns out that they are congruent via homothety.  To see this,
note that if~$X$ is a vector field on~$\bbC^2\simeq\bbP^2_0$ that 
generates dilation about a fixed point, then~$X$ lifts to a vector field~$Y$
on~$G_0$ that satisfies
$$
\Lie_Y\tau=\Lie_Y\phi=\Lie_Y\sigma=0,\qquad 
\Lie_Y\eta = c\,\eta,\quad \Lie_Y\omega = c\,\omega.
$$
for some nonzero (real) constant~$c$.  The vector field~$Y$ can then be
lifted to a vector field~$Z$ on~$G_0\times\bbC\times\bbC\times\bbR$ so
that it satisfies the same equations above as~$Y$ does but also satisfies
$$
\Lie_Z\abold = -c\,\abold,\quad \Lie_Z\sbold = -c\,\sbold,\quad 
\Lie_Z\pbold = -c\,\pbold.
$$
It then follows from the formulae for the generators of~${\Cal I}_3$ that
the flow of~$Z$ leaves~${\Cal I}_3$ invariant and therefore permutes the
leaves of~${\Cal I}_3$.  Since~$\Fbold$ is homogeneous of degree~$8$ in
$(\abold,\sbold,\pbold)$, it follows that~$\Lie_Z\Fbold = -8c\,\Fbold$.
In particular, the flow of~$Z$ acts as homothety on the level sets of~$\Fbold$.
Thus, any two solutions on which~$\Fbold$ is positive are 
congruent via homothety in~$\bbC^2\simeq\bbP^2_0$.

In this situation, it is therefore reasonable to restrict attention to the 
leaves that lie in the level set~$\Fbold=1$.  It is the solution corresponding 
to such a leaf that I am now going to describe.  Note that the isometry group
preserving such a solution is necessarily discrete.

\subsubhead 3.3.4. Local integration of the equations
\endsubsubhead
Suppose that one has a minimal Levi-flat hypersurface~$\Sigma\subset\bbC^2$
for which the bundle~$B_1^{**}$ satisfies equations~$(3)$.  Since~$s$ is
nonzero on~$B_1^{**}$, the structure equations show that there is a
submanifold~$B_2\subset B_1^{**}$ defined as the set on which~$s$ is 
real and positive and that~$B_2$ is a 6-fold cover of~$\Sigma^{**}$. From 
now on, all functions and forms are to be regarded as pulled 
back to~$B_2$. 

The reality of~$s$ and the structure equation
$$
s^{-1}\,ds = 6i\tau + p\,\eta + 3a\,\omega + 3 \bar a\,\bar\omega
$$
imply that~$\tau=0$.   Combining this with the structure equation
$d\eta = (a\,\omega{+}\bar a\,\bar\omega)\w\eta$ yields
$$
d\bigl( s^{-1/3}\,\eta\bigr) = 0.
$$
The structure equations also imply that the quantities
$$
\aligned
A &= {\ts{\frac19}}\,s^{2/3}
                \bigl(48\,| a|^2+12\,s^2+p^2\bigr),\\
B &= {\ts{\frac1{27}}} s\,
                \bigl(216\, a^2s + 216\,{\bar a}^2 s
                      +72\,| a|^2 p - 36\,s^2 p + p^3\bigr)\\ 
\endaligned
$$
introduced earlier satisfy equations of the form
$$
\aligned
dA &= -4\,B\phantom{^2}\,s^{-1/3}\,\eta \\
dB &= -6\,A^2\,s^{-1/3}\,\eta\\
\endaligned
$$
The assumption that the hypersurface~$\Sigma$
correspond to an ${\Cal I}_3$@-leaf on which~$\Fbold$ is identically equal
to~$1$ is equivalent to the equation~$A^3 - B^2 =1$, so there is a unique 
function~$\theta$ on~$B_2$ with values 
in the open interval~$(-\frac\pi2,\frac\pi2)$ for which
$$
A = \sec^{2/3}\theta > 0
\qquad\text{and}\qquad 
B = -\tan \theta.
$$
The above differential equations for~$A$ and~$B$ now imply that
$$
s^{-1/3}\,\eta = {\ts\frac16}\sec^{2/3}\theta\,d\theta.
$$
In particular, $d\theta$  never vanishes on~$B_2$ and is a nonzero multiple
of~$\eta$.

The structure equations now imply that
$$
d(s^{1/3}\,\omega) 
= \bigl(2\,s^{5/3}\,\bar\omega - {\ts{\frac13}}\,p\,s^{2/3}\,\omega\bigr)
 \w {\ts\frac16}\sec^{2/3}\theta\,d\theta.
$$
It follows that any point~$q\in B_2$ has a neighborhood~$U_0$ on which there 
exists a complex valued function~$z$, uniquely defined up to the addition 
of a (complex) function of~$t$, and a complex function~$L$, uniquely defined 
once~$z$ is chosen, so that
$$
s^{1/3}\,\omega = dz + {\ts\frac16}L\,\sec^{2/3}\theta\,d\theta.
$$
(Introducing such a coefficient in the~$L\,d\theta$ term simplifies 
later calculations.)  Because~$dz\w d\bar z \w d\theta 
= s^{2/3} \omega\w\bar\omega\w d\theta\not=0$, 
it follows that~$(z,\theta):U_0\to\bbC\times\bbR$ is a local diffeomorphism.  
By restricting to an appropriate neighborhood~$U_1\subset U_0$ of~$q$, I can
assume that~$(z,\theta):U_1\to\bbR\times\bbC$ defines a rectangular 
coordinate system (not necessarily centered on~$q$). 
Write~$z = x+i\,y$ where~$x$ and $y$ are real-valued.
Given the ambiguities in the choice of the coordinate system, 
partial differentiation with respect to~$z$ (or~$x$ or~$y$) is coordinate 
independent although partial differentiation with respect to~$\theta$ is not.

In these coordinates, the above structure equation 
for~$d\bigl( s^{1/3}\,\omega\bigr)$ now becomes
$$
dL\w \sec^{2/3}\theta\,d\theta
= \bigl(2\,s^{4/3}\,d\bar z-{\ts{\frac13}}\,p\,s^{1/3}\,dz\bigr)
          \w \sec^{2/3}\theta\,d\theta,
$$
so
$$
L_z = - {\ts{\frac13}}\,p\,s^{1/3}
\qquad\text{and}\qquad  
L_{\bar z} = 2\,s^{4/3}.
$$
Set $u=s^{1/3}$.  Then
$$
u^{-1}\,du = {\ts{\frac13}} s^{-1}\,ds 
    \equiv a\,\omega + \bar a\,\bar\omega
    \equiv u^{-1}\bigl( a\,dz + \bar a\, d{\bar z}\bigr) \bmod d\theta,
$$
so $a = u_z$.  The structure equation for~$da$ now gives
$$
da 
 \equiv (u^3 p - 3 a^2) u^{-1}\,dz + (|a|^2 -2u^6) u^{-1}\,d{\bar z} 
 \bmod d\theta,
$$
so it follows that $u_{zz} = a_z = \bigl(u^3 p -3{u_z}^2\bigr)u^{-1}$, which
can be written in the form
$$
(u^4)_{zz} = 4u^5\,p.
$$
Since~$p$ is real and since~$v_{zz} = {\ts{\frac14}}(v_{xx}-v_{yy}) 
+ {\frac{i}2}\,v_{xy}$ for any function~$v$ on~$U_1$, 
it follows that $(u^4)_{xy}=0$.  Consequently, there exist 
functions~$f$ and~$g$ defined on the rectangles~$(x,\theta)(U_1)$
and~$(y,\theta)(U_1)$ in~$\bbR^2$ so that
$$
u^4 = f(x,\theta) - g(y,\theta) > 0.
$$
These functions are unique up to the addition of a function of~$\theta$, i.e.,
one could replace~$\bigl(f(x,\theta),g(y,\theta)\bigr)$ 
by~$\bigl(f(x,\theta){+}h(\theta),g(y,\theta){+}h(\theta)\bigr)$ 
for some~$h$ defined on the interval~$\theta(U_1)$, 
but this is the only ambiguity in the choice of these two functions.

Now, the equation for~$da$ also implies the equation~$u_{z\bar z} = a_{\bar z}
= (|u_z|^2 - 2\,u^6) u^{-1}$, which can be written in the form
$$
(u^4)_{z\bar z} = u^{-4}\,\bigl|(u^4)_z\bigr|^2 - 8\,(u^4)^2.
$$
Using the expression already found for~$u^4$ plus the formulae
$v_{z\bar z} = {\ts{\frac14}}(v_{xx}+v_{yy})$ and $|v_z|^2 
= {\ts{\frac14}}\bigl({v_x}^2+{v_y}^2\bigr)$, this equation can be written
in the form
$$
f_{xx}(x,\theta) - g_{yy}(y,\theta) 
     = \frac{{f_x(x,\theta)}^2+{g_y(y,\theta)}^2}{f(x,\theta) - g(y,\theta)}
           -32\bigl(f(x,\theta) - g(y,\theta)\bigr)^2.
$$

Now,  setting~$v = f-g$, this can be written in the form
$$
v_{xx} + v_{yy} = \frac{{v_x}^2+{v_y}^2}{v}
           -32\,v^2
$$
and rearranged to give
$$
\left(\frac{v_x}{v}\right)_x 
= \frac{{v_y}^2}{v^2} + \frac{v_{yy}}{v} - 32\,v.
$$
Since~$v_{xy}=0$, both~$v_y$ and $v_{yy}$ are constant in~$x$.  
Thus, multiplying this equation by~$2v_x/v$ and 
integrating with respect to~$x$ yields
$$
\left(\frac{v_x}{v}\right)^2 
= C(y,\theta) - \frac{{v_y}^2}{v^2} - \frac{2v_{yy}}{v} - 64\,v.
$$
for some function~$C$ on~$(y,\theta)(U_1)$.  
Now, multiplying by~$v^2$
and substituting~$v = f - g$, this can be written in the form
$$
{f_x(x,\theta)}^2  = 2\,a_0(y,\theta) 
                        +12\,a_1(y,\theta)\,f(x,\theta)
                        +48\,a_2(y,\theta)\,{f(x,\theta)}^2 
                        -64\,{f(x,\theta)}^3
$$
for some functions~$a_0$, $a_1$, and $a_2$ on~$(y,\theta)(U_1)$. 
(The choice of numerical coefficients is cosmetic.)  

Now, if the functions~$a_i$ really did depend on~$y$, differentiating
this equation with respect to~$y$ would then force~$f(x,\theta)$ 
to be constant in~$x$, making~$f_x$ vanish identically.  This would,
in turn, imply that~$a = u_z = -{\ts\frac12}i\,u_y$ is purely imaginary,
so that the quantity~$a^2s$ would be real.  However, going back to the 
analysis in~\S2.2.2, this can only happen for solutions of type~$(2)$.
Since the goal of this section is analyse the solutions of type~$(3)$ 
that have not already been accounted for by those of type~$(1)$ or~$(2)$,
this case can therefore be set aside.

Thus, $f$ satisfies an equation of the form
$$
{f_x(x,\theta)}^2  = 2\,a_0(\theta) 
                        +12\,a_1(\theta)\,f(x,\theta)
                        +48\,a_2(\theta)\,{f(x,\theta)}^2 
                        -64\,{f(x,\theta)}^3
$$
for some functions~$a_0$, $a_1$, and $a_2$ on~$\theta(U_1)$.  A similar
analysis shows that there are functions~$b_0$, $b_1$, and $b_2$ 
on~$\theta(U_1)$ for which
$$
{g_y(y,\theta)}^2  = -2\,b_0(\theta) 
                        -12\,b_1(\theta)\,g(y,\theta)
                        -48\,b_2(\theta)\,{g(y,\theta)}^2 
                        +64\,{g(y,\theta)}^3.
$$
Moreover, substituting these relations and their derivatives 
back into the original equation for~$v$, it follows 
that~$b_0=a_0$, $b_1=a_1$, and $b_2=a_2$.  Thus,
$$
\alignedat5
{f_x(x,\theta)}^2 &\ =\ &\ \phantom{-}2\,a_0(\theta) 
                        &\ +\ & 12\,a_1(\theta)\,f(x,\theta)
                        &\ +\ & 48\,a_2(\theta)\,{f(x,\theta)}^2 
                        &\ -\ & 64\,{f(x,\theta)}^3,& \\
{g_y(y,\theta)}^2 &\ =\ &\           -2\,a_0(\theta) 
                        &\ -\ & 12\,a_1(\theta)\,g(y,\theta)
                        &\ -\ & 48\,a_2(\theta)\,{g(y,\theta)}^2 
                        &\ +\ & 64\,{g(y,\theta)}^3.& \\
\endalignedat
$$
By replacing~$\bigl(f(x,\theta),g(y,\theta)\bigr)$ 
with 
$\bigl(f(x,\theta){-}\frac14\,a_2(\theta),
g(y,\theta){-}\frac14\,a_2(\theta)\bigr)$, 
it can be arranged that~$a_2\equiv0$.   This removes the ambiguity in
the choice of~$f$ and~$g$.  

At this point, $f$ and $g$ satisfy the equations
$$
\alignedat4
{f_x(x,\theta)}^2 
&\ =\ &\ \phantom{-}2\,a_0(\theta) &\ +\ & 12\,a_1(\theta)\,f(x,\theta)
                      &\ -\ & 64\,{f(x,\theta)}^3& \\
{g_y(y,\theta)}^2 
&\ =\ &\           -2\,a_0(\theta) &\ -\ & 12\,a_1(\theta)\,g(y,\theta)
                      &\ +\ & 64\,{g(y,\theta)}^3& \\
\endalignedat
$$
as well as equations
$$
\alignedat3
f_{xx}(x,\theta) 
&\ =\ &\phantom{-}6\,a_1(\theta) &\ -\ & 96\,{f(x,\theta)}^2\,,& \\
g_{yy}(y,\theta) 
&\ =\ &         - 6\,a_1(\theta) &\ +\ & 96\,{g(y,\theta)}^2\,.& \\
\endalignedat
$$
This information can now be substituted back 
into the previous formulae, yielding
$$
\aligned
s &= \bigl(f(x,\theta)-g(y,\theta)\bigr)^{3/4}\\
a &={\ts\frac18}\,\bigl(f(x,\theta)-g(y,\theta)\bigr)^{-3/4}
                       \bigl(f_x(x,\theta)+i\,g_y(y,\theta)\bigr)\\
p &= -6\,\,\bigl(f(x,\theta)-g(y,\theta)\bigr)^{-1/4}
      \bigl(f(x,\theta)+g(y,\theta)\bigr).\\
\endaligned
$$
Using these formulae, the definitions of~$A$ and~$B$, and
the equations satisfied by~$f$ and $g$, 
it now follows that
$$
a_1(\theta) = A = \sec^{2/3}\theta 
\qquad\text{and}\qquad 
a_0(\theta) = B = -\tan\theta.
$$

The previous formula for~$dL$ now simplifies to
$$
dL\equiv 4\,f(x,\theta)\,dx + 4\,g(y,\theta)\,dy \bmod d\theta,
$$
so that~$L = F(x,\theta)+i\,G(y,\theta)$ for functions~$F$ and~$G$ 
satisfying~$F_x = 4 f$ and $G_y = 4 g$.  All this information
combines to yield the formulae 
$$
\aligned
\tau &= 0\\
\eta &= {\ts\frac16}\bigl(f(x,\theta)-g(y,\theta)\bigr)^{1/4}\,
                \sec^{2/3}\theta\,d\theta\\
\omega &= \bigl(f(x,\theta)-g(y,\theta)\bigr)^{-1/4}
         \,\bigl(dz + {\ts\frac16}\sec^{2/3}\theta\,
                 \bigl(F(x,\theta)+i\,G(y,\theta)\bigr)\,d\theta\bigr).\\
\endaligned
$$

Now, the cubic polynomial
$$
p(\lambda,\theta) 
= -2\,\tan\theta +12\,\sec^{2/3}\theta\,\lambda - 64\,\lambda^3
$$
has three real, distinct roots in~$\lambda$.
In fact, defining
$$
\aligned
r_1(\theta) &= {\ts\frac12}\,\sin\bigl({\ts\frac13}\,\theta
                                         - {\ts\frac23}\pi\bigr)
                                                  \,\sec^{1/3}\theta,\\
r_2(\theta) &= {\ts\frac12}\,\sin\bigl({\ts\frac13}\,\theta
                              \phantom{{}+ {\ts\frac03}\pi}\bigr)
                                                  \,\sec^{1/3}\theta,\\
r_3(\theta) &= {\ts\frac12}\,\sin\bigl({\ts\frac13}\,\theta 
                                         + {\ts\frac23}\pi\bigr)
                                                  \,\sec^{1/3}\theta,\\
\endaligned
$$
one has $r_1(\theta)<r_2(\theta)<r_3(\theta)$ 
when~$-\frac\pi2<\theta<\frac\pi2$ and
$$
p(\lambda,\theta) = -64\,\bigl(\lambda - r_1\bigl(\theta\bigr)\bigr)
                         \bigl(\lambda - r_2\bigl(\theta\bigr)\bigr)
                         \bigl(\lambda - r_3\bigl(\theta\bigr)\bigr).
$$

Now, the differential equations on~$f(x,\theta)$ and $g(y,\theta)$ coupled 
with the inequality~$g(x,\theta)<f(y,\theta)$ imply the inequalities
$$
r_1(\theta)<g(y,\theta) < r_2(\theta) < f(x,\theta) < r_3(\theta).
$$
Moreover the differential equation for~$f$ (resp.~$g$) can now be used to
extend its range of definition from~$(x,\theta)(U_1)$ (resp.~$(y,\theta)(U_1)$)
to all of~$\bbR\times(\frac\pi2,\frac\pi2)$.  The extended functions satisfy
$$
r_1(\theta) \le g(y,\theta) \le  r_2(\theta) \le  f(x,\theta) \le r_3(\theta)
$$
and the periodicity relations
$$
\aligned
f\bigl(x+2\rho_+(\theta),\theta\bigr) &= f(x,\theta)\\
g\bigl(y+2\rho_-(\theta),\theta\bigr) &= g(y,\theta)\\
\endaligned
$$
where the functions~$\rho_\pm$ are defined by the elliptic integrals 
$$
\aligned
\rho_+(\theta) &= \frac18\,\int_{r_2(\theta)}^{r_3(\theta)}
\frac{da}{\sqrt{\bigl(r_3(\theta)-a\bigr)
                 \bigl(a-r_2(\theta)\bigr)\bigl(a-r_1(\theta)\bigr)}}\,,\\
\rho_-(\theta) &= \frac18\,\int_{r_1(\theta)}^{r_2(\theta)}
\frac{da}{\sqrt{\bigl(r_3(\theta)-a\bigr)
                 \bigl(r_2(\theta)-a\bigr)\bigl(a-r_1(\theta)\bigr)}}\,.\\
\endaligned
$$
Note, by the way, that~$\rho_+(-\theta) = \rho_-(\theta)>0$ for~$\theta$
in~$(-\frac\pi2,\frac\pi2)$.  

Using these extended functions, I can now modify $x$ and $y$ by adding 
functions of~$\theta$ so as to arrange that
$$
g(0,\theta) = r_2(\theta) = f(0,\theta).
$$
This makes the coordinates~$(x,y,\theta)$ unique up to replacement by
coordinates of the form
$$
(x^*,y^*,\theta^*) 
= \bigl(\,x+2m\,\rho_+(\theta),\, y+2n\,\rho_-(\theta),\,\theta\,\bigr)
$$
for some integers~$m$ and $n$.  These formulae will be important in
the discussion of discrete symmetries that will be undertaken below.

The functions~$f$ and $g$ are now uniquely defined on the entire
strip~$\bbR\times(-\frac\pi2,\frac\pi2)$ by the requirement that they 
satisfy the second order equations with initial conditions
$$
\alignedat3
f_{xx}(x,\theta) &= \phantom{-}6\,\sec^{2/3}\theta - 96\,f(x,\theta)^2\,,
   &\quad f(0,\theta) &= r_2(\theta), &\quad f_x(0,\theta) &= 0,\\
g_{yy}(y,\theta) &=           -6\,\sec^{2/3}\theta + 96\,g(y,\theta)^2\,,
   &\quad g(0,\theta) &= r_2(\theta), &\quad g_y(0,\theta) &= 0.\\
\endalignedat
$$
Then~$u(x,y,\theta)^4 = f(x,\theta)-g(y,\theta)\ge0$ 
is doubly periodic on~$\bbR\times\bbR\times(-\frac\pi2,\frac\pi2)$ in the 
obvious sense and is strictly positive except along the curves~$C_{m,n}$ of 
the form~$(x,y,\theta) = \bigl(2m\rho_+(\theta),2n\rho_-(\theta),\theta\bigr)$
for any integers~$m$ and~$n$.  The vanishing near these lines is very
simple:  Along~$C_{0,0}$, i.e., the line~$(x,y,\theta)=(0,0,\theta)$, there
are convergent Taylor expansions
$$
f(x,\theta) = r_2(\theta) + \sum_{k=1}^\infty c_k(\theta)\,x^{2k}\,,
\qquad\qquad
g(y,\theta) = r_2(\theta) + \sum_{k=1}^\infty (-1)^k\,c_k(\theta)\,y^{2k}\,,
$$
implying that there is a smooth function~$\tilde u$ 
on~$\bbR\times\bbR\times(-\frac\pi2,\frac\pi2)$ 
satisfying~$\tilde u(0,0,\theta) 
= c_1(\theta) = 3\sec^{2/3}\theta\,\bigl(1-4\sin^2(\frac13\theta)\bigr)>0$ 
for which
$$
u(x,y,\theta)^4 = (x^2+y^2)\,{\tilde u}(x,y,\theta).
$$
By the periodicity relations, the description of the
vanishing of~$u$ near the other curves~$C_{m,n}$ follows from this one.

Now, examining the coefficient of~$d\theta$ in the formula for~$ds$ yields 
the relation
$$
  g_y\,G - 6\,\cos^{2/3}\theta\,g_\theta - 8\,g^2 
= f_x\,F - 6\,\cos^{2/3}\theta\,f_\theta - 8\,f^2.
$$
The left hand side of this relation is independent of~$x$ while the
right hand side is indepdendent of~$y$, so that each side is
a function of~$\theta$ only.  Evaluating either side at~$x=y=0$ then yields
$$
  g_y\,G - 6\,\cos^{2/3}\theta\,g_\theta - 8\,g^2 
= f_x\,F - 6\,\cos^{2/3}\theta\,f_\theta - 8\,f^2
= -6\,\cos^{2/3}\theta\,r'_2(\theta) - 8\,r_2(\theta)^2.
$$
Of course, this allows one to solve for~$F$ and $G$ away from the places
where~$f_x$ and $g_y$ vanish, yielding formulae of the form
$$
\aligned
F &= \bigl[6\,\cos^{2/3}\theta\,\,(\,f_\theta-r'_2(\theta)\,) 
                    +8\,(f^2-r_2(\theta)^2)\bigr]/f_x\\
G &= \bigl[6\,\cos^{2/3}\theta\,\,(\,g_\theta-r'_2(\theta)\,) 
                    +8\,(g^2-r_2(\theta)^2)\bigr]/g_y\\
\endaligned
$$
Since~$f_x(x,\theta)=0$ if and only if~$x$ is an integer 
multiple of~$\rho_+(\theta)$ and~$g_y(y,\theta)=0$ 
if and only if~$y$ is an integer multiple of~$\rho_-(\theta)$,
this gives integration-free formulae for~$F$ and $G$ that 
are valid over a dense open set.  Moreover, differentiating 
the relations above with respect to~$x$ or~$y$
and using the identities~$F_x = 4f$ and~$G_y=4g$ yields
$$
  g_{yy}\,G - 6\,\cos^{2/3}\theta\,g_{y\theta} - 12\,g g_y 
= f_{xx}\,F - 6\,\cos^{2/3}\theta\,f_{x\theta} - 12\,f f_y
= 0.
$$
Since~$f_x$ and $f_{xx}$ do not vanish simultaneously, 
and since~$g_y$ and $g_{yy}$ do not vanish simultaneously, 
these relations together with the relations above yield 
explicit smooth formulae for~$F$ and~$G$ over all 
of~$\bbR\times(-\frac\pi2,\frac\pi2)$.  In particular,
these formulae imply that~$F(0,\theta)=G(0,theta)=0$, so
that~$F$ and $G$ can also be described by
$$
F(x,\theta) = 4\int_0^xf(\xi,\theta)\,d\xi,\qquad\qquad
G(y,\theta) = 4\int_0^yf(\xi,\theta)\,d\xi.
$$

The integration-free formulae yield pseudo-periodicity relations 
for~$F$ and~$G$:  Differentiating
$$
f\bigl(x+2\rho_+(\theta),\theta\bigr) = f(x,\theta)
$$
with respect to~$\theta$ shows that~$f_\theta$ satisfies the 
pseudo-periodicity relation
$$
f_\theta\bigl(x+2\rho_+(\theta),\theta\bigr) - f_\theta(x,\theta)
= -2f_x(x,\theta)\,\rho'_+(\theta).
$$
Consequently, $F$ satisfies the pseudo-periodicity relation
$$
F\bigl(x+2\rho_+(\theta),\theta\bigr) - F(x,\theta)
= -12\,\rho'_+(\theta)\,\cos^{2/3}\theta\,.
$$
Similarly,
$$
G\bigl(y+2\rho_-(\theta),\theta\bigr) - G(y,\theta)
= -12\,\rho'_-(\theta)\,\cos^{2/3}\theta\,.
$$


At this point, all the structure equations in~$(3)$ are identities.

\subsubhead 3.3.5. Global structure of the solution
\endsubsubhead
The local information derived in the previous subsubsection can
now be used to give a global description of the corresponding minimal
Levi-flat hypersurface in~$\bbC^2$.  To begin, define~$r_i$ for~$i=1$, $2$, 
and $3$ and $\rho_\pm$ as functions on~$(-\frac\pi2,\frac\pi2)$ by the 
already listed formulae.  Then, define functions~$f$ and $g$ 
on~$\bbR\times(-\frac\pi2,\frac\pi2)$ by the differential equations
with initial conditions:
$$
\alignedat3
f_{xx}(x,\theta) &= \phantom{-}6\,\sec^{2/3}\theta - 96\,f(x,\theta)^2\,,
   &\quad f(0,\theta) &= r_2(\theta), &\quad f_x(0,\theta) &= 0,\\
g_{yy}(y,\theta) &=           -6\,\sec^{2/3}\theta + 96\,g(y,\theta)^2\,,
   &\quad g(0,\theta) &= r_2(\theta), &\quad g_y(0,\theta) &= 0.\\
\endalignedat
$$
Note that~$f$ is even and periodic of period~$2\rho_+(\theta)$ in its
first argument while $g$ is even and periodic
of period~$2\rho_-(\theta)$ in its first argument.  
Moreover, these functions automatically satisfy the 
first order equations
$$
\alignedat4
{f_x(x,\theta)}^2 &\ =\ &\ \phantom{-}2\,\tan\theta 
                     &\ +\ & 12\,\sec^{2/3}\theta\,f(x,\theta)
                      &\ -\ & 64\,{f(x,\theta)}^3,& \\
{g_y(y,\theta)}^2 &\ =\ &\           -2\,\tan\theta 
                     &\ -\ & 12\,\sec^{2/3}\theta\,g(y,\theta)
                      &\ +\ & 64\,{g(y,\theta)}^3.& \\
\endalignedat
$$
Define~$F$ and~$G$ on the same domain by
$$
F(x,\theta) = \int_0^x 4\,f(\xi,\theta)\,d\xi,
\qquad\qquad
G(y,\theta) = \int_0^y 4\,g(\xi,\theta)\,d\xi.
$$

Let~$D =  \bbR\times\bbR\times(-\frac\pi2,\frac\pi2)$ 
and let~$D^*\subset D$ be the complement of the curves
$$
C_{m,n} = \left\{\,\bigl(2m\rho_+(\theta),2n\rho_-(\theta),\theta\bigr)
               \ \vrule\ {\ts\theta\in(-\frac\pi2,\frac\pi2)}\ \right\}.
$$
Finally, define functions and 1-forms on~$D^*$ by the formulae
$$
\aligned
s &= \bigl(f(x,\theta)-g(y,\theta)\bigr)^{3/4},\\
a &={\ts\frac18}\,\bigl(f(x,\theta)-g(y,\theta)\bigr)^{-3/4}
                    \bigl(f_x(x,\theta)+i\,g_y(y,\theta)\bigr),\\
p &= -6\,\,\bigl(f(x,\theta)-g(y,\theta)\bigr)^{-1/4}
                     \bigl(f(x,\theta)+g(y,\theta)\bigr),\\
\eta &= {\ts\frac16}\,\bigl(f(x,\theta)-g(y,\theta)\bigr)^{1/4}\,
                   \sec^{2/3}\theta\,d\theta,\\
\omega &= \bigl(f(x,\theta)-g(y,\theta)\bigr)^{-1/4}
               \,\bigl(dz  +{\ts\frac16}\sec^{2/3}\theta\,
          \bigl(F(x,\theta)+i\,G(y,\theta)\bigr)\,d\theta\bigr),\\
\tau &=0.\\
\endaligned
$$
Then the structure equations~$(3)$ are satisfied on~$D^*$. In particular,
setting $\sigma = -2\,\bar a\,\eta - 2\,s\,\bar\omega$ and
$\phi = -i\,a\,\omega + i\,\bar a\,\bar\omega$, 
the $\eug_0$@-valued 1-form
$$
\gamma = 
\pmatrix
      0    &     0   &      0      \\
      \eta & i\,\phi & -\bar\sigma \\
    \omega & \sigma  & -i\,\phi    \\
\endpmatrix
$$
satisfies~$d\gamma = -\gamma\w\gamma$.

By the usual moving frame argument~\cite{Gr}, it follows that,
if~$U\subset D^*$ is any simply connected domain in~$D^*$, 
then there is a map~$\gbold:U\to G_0$, unique up to left translation 
by a constant, so that~$\gbold^{-1}\,d\gbold = {\gamma\vrule_{U}}$.  
The projection~$\gbold K:U\to G_0/K=\bbC^2$ is then an immersion of~$U$ 
into~$\bbC^2$ as a minimal Levi-flat hypersurface of type~$(3)$.  However, 
this argument does not provide a description of the topology 
or global properties of the solution.  It is to this description that 
I now turn.

The group~$\bbZ^2$ acts on~$D$ preserving~$D^*$ via the maps
$$
\Phi_{m,n}(x,y,\theta) = 
\bigl( x + 2m\,\rho_+(\theta),\,y + 2n\,\rho_-(\theta),\,\theta\bigr).
$$
Denote the~$\bbZ^2$-orbit 
of~$(x,y,\theta)$ by~$[x,y,\theta]\in N$.  
The periodicity relations on~$f$ and $g$ combined with the 
pseudo-periodicity relations on~$F$ and~$G$ imply~$\Phi_{m,n}^*\gamma 
= \gamma$.  (In fact, all the quantitites~$s$, $a$, $p$, $\eta$, $\omega$, 
and~$\tau$ $(=0)$ are invariant under this~$\bbZ^2$-action.)
Thus~$\gamma$ is well-defined on the quotient space~$N^* =
D^*/\bbZ^2$, which is diffeomorphic to a punctured torus cross an
open interval. 

On~$N^*\times G_0$, thought of as a trivialized principal left
$G_0$-bundle over~$N^*$, consider the $\eug_0$-valued 
connection $1$-form
$$
\psi = dg\,g^{-1} -g\,\gamma\,g^{-1}
     = g\bigl(\,g^{-1}\,dg - \gamma\bigr)\,g^{-1}.
$$  
Since~$d\gamma = -\gamma\w\gamma$, it follows that~$d\psi = \psi\w\psi$,
i.e., that~$\psi$ is flat. 
Consequently, $N^*\times G_0$ is foliated by $\psi$-leaves, each of
which is a smooth submanifold~$L\subset N^*\times G_0$ such that projection
onto the first factor is a covering map and such that any two leaves differ
by left action in the $G_0$-factor by a constant element of~$G_0$.  
For any such leaf~$L$, we can regard the functions~$s$, $a$, $p$ and 
$1$@-forms $\eta$ and $\omega$ as being well defined on~$L$ via pullback
from the projection~$L\to N^*$.  

The map~$(g,a,s,p):L\to M^{13}=G_0\times\bbC
\times\bbC\times\bbR$ then immerses~$L$ as an $\cI_3$@-leaf lying 
in the locus~$\Fbold=1$.  By the construction of~$\gamma$ and the development 
that led up to it, the image of~$L$ is a complete $\cI_3$@-leaf.  
Thus, the topology of the leaves will be known once the 
covering map~$L\to D^*$  and the projection~$g:L\to G_0$ are understood.

The projection~$g:L\to G_0$ is simply a diffeomorphism.  This follows
because, on~$L$, the $\eug_0$-valued 1-form~$g^{-1}\,dg$ is simply~$\gamma$,
which determines the forms~$\omega$ and~$\eta$ and the functions~$s$, $a$,
and~$p$.  The construction of the coordinate system~$(x,y,\theta)$
from~$(\eta,\omega,s,a,p)$ shows that this suffices to recover
the map~$(x,y,\theta):L\to D^*$ up to the action of~$\bbZ^2$, which is
the same as recovering~$[x,y,\theta]:L\to N^*$ and hence the full embedding
of~$L$ into~$N^*\times G_0$.   In particular, this implies that~$(g,a,s,p)$
is an embedding. 

Now, a leaf~$L$ is just the holonomy bundle of~$\psi$ 
through each of its points.  For the sake of concreteness, 
choose~$n_0 = [\rho_0,\rho_0,0]\in N^*$ as basepoint, 
where~$\rho_0 = \rho_+(0)=\rho_-(0)$ and let~$L\subset N^*\times G_0$ 
be the leaf of~$\psi$ that passes 
through~$\bigl(n_0,\text{I}_3\bigr)$.  The 
intersection~$L\cap\bigl(\{n_0\}\times G_0\bigr)$ is then of the form
$\{n_0\}\times \Gamma$ where~$\Gamma\subset G_0$ 
is the holonomy subgroup of~$\psi$ and this is what must be computed.
The calculations below will actually determine the $\psi$-monodromy 
homomorphism~$\pi_1\bigl(N^*,n_0)\to G_0$, whose image is~$\Gamma$.

Since~$N^*$ is an interval cross a punctured torus, $\pi_1\bigl(N^*,n_0)$
is generated by the loops~$X:\bigl[0,2\rho_0\bigr]\to N^*$ 
and $Y:\bigl[0,2\rho_0\bigr]\to N^*$ defined by
$$
X(x) = \bigl[x{+}\rho_0,\,\rho_0,\,0\bigr],\qquad\qquad\quad
Y(y) = \bigl[\rho_0,\,y{+}\rho_0,\,0\bigr].
$$

To compute the $\psi$-monodromy around these two loops, information
about the behavior of the functions~$f$ and $g$ when $\theta=0$ will be 
used. To begin, note that~$r_1(0) = -\sqrt3/4$, $r_2(0)=0$, and 
$r_3(0)=\sqrt3/4$ and observe that, by the symmetry properties 
of~$f$ and~$g$, there is a $2\rho_0$@-periodic function~$v$ 
on~$\bbR$ that satisfies
$$
v(t) = f(t{+}\rho_0,0)+\sqrt3/4 = \sqrt3/4 - g(t{+}\rho_0,0)
$$
for all~$t$.  In fact, $v$ is defined by the conditions that it
satisfy both the initial condition~$v(0) = \sqrt3/2$ and the 
Weierstra\ss-type differential equation
$$
\bigl(v'(t)\bigr)^2 
  = 64\,v(t)\bigl(\sqrt3/2-v(t)\bigr)
                  \bigl(v(t)-\sqrt3/4\bigr).
$$
Note that~$v$ is positive, satisfying~$\sqrt3/4\le v(t)
\le \sqrt3/2$, and that~$v$ is an even function on~$\bbR$. 
In particular, satisfies~$v(2\rho_0-t) = v(t)$, 
a fact that will be used below.


Now, from the definition of~$X$ it follows that
$$
X^*(\gamma) = 
\pmatrix
      0    &     0   &      0      \\
      0 & 0 & 2\bigl(v(x)\bigr)^{1/2}\,dx \\
    \bigl(v(x)\bigr)^{-1/4}\,dx & -2\bigl(v(x)\bigr)^{1/2}\,dx & 0 \\
\endpmatrix,
$$
Consider the~$g_X:\bigl[0,2\rho_0\bigr]\to G_0$ that
satisfies~${g_X}^{-1}\,dg_X = X^*\gamma$ and $g_X(0)=\text{I}_3$.  
Because~$X^*\gamma$ takes values in~$\eug_0\cap\eusl(3,\bbR)$, the
map~$g_X$ has values in~$G_0\cap\SL(3,\bbR)$ and so can be
written in the form
$$
g_X(x) = \pmatrix 1 & 0 & 0\\ 
                  u_1(x) & \phantom{-}\cos\varphi(x) & \sin\varphi(x)\\
                  u_2(x) & -\sin\varphi(x) & \cos\varphi(x)\\
          \endpmatrix,
$$
where the functions~$u_1$, $u_2$, and $\varphi$ on~$[0,2\rho_0]$ are defined
by the ODE system
$$
\alignedat2
u_1'(x) &= \sin\varphi(x)\, \bigl(v(x)\bigr)^{-1/4},&
          \qquad\qquad u_1(0) & = 0,\\
u_2'(x) &= \cos\varphi(x)\, \bigl(v(x)\bigr)^{-1/4},&
          \qquad\qquad u_2(0) & = 0,\\
\varphi'(x) 
  &= \phantom{\sin\varphi(x)}\llap{2}\,
      \bigl(v(x)\bigr)^{1/2\phantom{-}},&
          \qquad\qquad \varphi(0) & = 0.\\
\endalignedat
$$

The ODE that~$v$ satisfies suggests a change 
of variables eliminating the explicit $x$-dependence,
yielding
$$
\align
\varphi\bigl(2\rho_0\bigr) &=
\int_0^{2\rho_0}2\bigl(v(x)\bigr)^{1/2}\,dx  \\
&= 4\cdot\frac18\int_{\sqrt3/4}^{\sqrt3/2} 
\frac{v^{1/2}\,dv}
  {\sqrt{v\bigl(\sqrt3/2-v\bigr)\bigl(v-\sqrt3/4)\bigr)}}\\
&=\frac\pi2.\\
\endalign
$$
Thus~$\varphi$ defines a 
diffeomorphism~$\varphi:\bigl[0,2\rho_0\bigr]\to[0,\frac\pi2]$ that, 
because of the symmetries of~$v$, 
has the symmetry~$\varphi\bigl(2\rho_0{-}x\bigr) = \frac\pi2-\varphi(x)$.  
In turn, this implies that~$u_i'(x)>0$ for all~$x\in\bigl(0,2\rho_0\bigr)$
and, by a straightforward change of variables, 
that~$u_1\bigl(2\rho_0\bigr)=u_2\bigl(2\rho_0\bigr) = r$ for some%
\footnote{For the curious: Numerical calculation yields
$\rho_0\approx 0.498083225$ and $r\approx .565201447$.} 
$r>0$.

This implies that~$g_X\bigl(2\rho_0\bigr)=h_X$ where
$$
h_X = \pmatrix1 & 0 & 0\\ r & 0 & 1\\ r &-1&0\\ \endpmatrix.
$$
This~$h_X$ represents the holonomy of~$\psi$ around the loop~$X$.  
(Note that it is possible to compute the map~$g_X$ and hence the
holonomy~$h_X$ by quadratures in this manner because~$X^*\gamma$ 
takes values in a solvable subalgebra of~$\eug_0$.)

A similar argument for~$Y$ gives
$$
Y^*(\gamma) = 
\pmatrix
      0    &     0   &      0      \\
      0 & 0 & -2i\,\bigl(v(y)\bigr)^{1/2}\,dy \\
    \bigl(v(y)\bigr)^{-1/4}\,dy & -2i\,\bigl(v(y)\bigr)^{1/2}\,dy & 0 \\
\endpmatrix,
$$
Carrying out the same sort of analysis as was applied to~$X$ leads to the
conclusion that if~$g_Y:\bigl[0,2\rho_0\bigr]\to G_0$ is the map that
satisfies~${g_Y}^{-1}\,dg_Y = Y^*\gamma$ and $g_Y(0)=\text{I}_3$, 
then~$g_Y\bigl(2\rho_0\bigr)=h_Y$ where
$$
h_Y = \pmatrix1 & 0 & 0\\ -ir & 0 & -i \\ r &-i&0\\ \endpmatrix.
$$
Thus~$h_Y$ represents the holonomy of~$\psi$
around the loop~$Y$.

Now, setting
$$
\vbold = \pmatrix 0\\ r\endpmatrix  \not=0,
$$
and
$$
\bone = \pmatrix 1&0\\ 0&1\endpmatrix,
\quad
\ibold = \pmatrix \phantom{-}0&1\\ -1&0\endpmatrix,
\quad 
\jbold = \pmatrix \phantom{-}0&-i\\ -i&\phantom{-}0\endpmatrix,
\quad
\kbold = \ibold\jbold 
       = \pmatrix -i&0\\ 0&i\endpmatrix,
$$
it follows that
$$
h_X = \pmatrix 1& 0\\ (\bone+\ibold)\vbold&\ibold\endpmatrix,
\qquad\qquad
h_Y = \pmatrix 1& 0\\ (\bone+\jbold)\vbold&\jbold\endpmatrix.
$$

Noting that~$\ibold^2=\jbold^2=-\bone$ while~$\kbold = \ibold\jbold 
= -\jbold\ibold$, it is evident that~${h_X}^4 = {h_Y}^4 = \text{I}_3$
and that any iterated product of the matrices~$h_X$ and~$h_Y$ is of
the form
$$
h = \pmatrix 1& 0\\ 
    (a_0\bone{+}a_1\ibold{+}a_2\jbold{+}a_3\kbold)\vbold&\qbold\endpmatrix
$$
where~$\qbold$ lies in~$\{\pm\bone,\pm\ibold,\pm\jbold,\pm\kbold\}$ and
the~$a_i$ are integers whose sum is even.  In particular, the 
subgroup~$\Gamma\subset G_0$ generated by~$h_X$ and $h_Y$ is discrete.
Moreover, the 
homomorphism~$\Gamma\to \{\pm\bone,\pm\ibold,\pm\jbold,\pm\kbold\}$
defined by~$h\mapsto \qbold$ in the above notation is surjective.
It is not difficult to establish that the kernel~$\hat\Lambda$ of this 
homomorphism consists exactly of the matrices of the form
$$
\pmatrix 1& 0\\ 
    2(a_0\bone{+}a_1\ibold{+}a_2\jbold{+}a_3\kbold)\vbold&\bone\endpmatrix
$$
where the~$a_i$ are integers whose sum is even.  Since~$\vbold\not=0$, 
the set~$\Lambda\subset\bbC^2$ consisting 
of the vectors~$2(a_0\bone{+}a_1\ibold{+}a_2\jbold{+}a_3\kbold)\vbold$ where 
the~$a_i$ are integers whose sum is even is a lattice in~$\bbC^2$, i.e.,
a discrete abelian subgroup of rank~$4$.  Up to rotation and scaling, 
$\Lambda$ is a lattice of type~$\text{F}_4$.  In what follows, it will
be useful to identify~$\Lambda$ with~$\hat\Lambda\subset G_0$ 
via the identification
$$
2(a_0\bone{+}a_1\ibold{+}a_2\jbold{+}a_3\kbold)\vbold
\longmapsto
\pmatrix 1& 0\\ 
    2(a_0\bone{+}a_1\ibold{+}a_2\jbold{+}a_3\kbold)\vbold&\bone\endpmatrix,
$$
so I will do this henceforth without explicit comment.

Let~$\hat K\subset\pi_1(N^*,n_0)$ denote the normal subgroup of 
index~$8$ consisting of those homotopy classes of loops whose 
$\psi$@-holonomy  lies in~$\hat\Lambda$ and let~$\hat N^*\to N^*$ denote 
the $8$@-fold covering space corresponding to~$\hat K$.  I am going 
to show that there is a way of `completing'~$\hat N^*$ in a natural way 
so that each of the complex leaves of~$\hat N^*$ (i.e., the leaves 
of~$\eta=0$) is realized as a compact Riemann surface of genus~$3$ 
punctured at four points. I will then examine to what extent the functions 
and forms~$s$, $a$, $p$, $\eta$, and~$\omega$ extend smoothly across these 
punctures.  

Ultimately, the goal is to show that~$\bbC^2/\Lambda$ contains a minimal
Levi-flat hypersurface whose complex leaves are (compact) Riemann surfaces
of genus~$3$. 

Let~$\tilde N$ be the quotient of~$D$ by the action of~$(2\bbZ)^2$, 
i.e., the index~$4$ subgroup of~$\bbZ^2$ generated by the 
transformations~$\Phi_{2m,2n}$, and let~$\tilde N^*\subset\tilde N$
be the image of~$D^*\subset D$ under this quotient action.
Let~$\la x,y,\theta\ra\in \tilde N$ denote the 
equivalence class of~$(x,y,\theta)\in D$ under the action of~$(2\bbZ)^2$.
Any product of a finite sequence drawn from~$\{h_X,h_Y\}$
that contains an odd number of copies of either~$h_X$ or $h_Y$ will be
an~$h\in \Gamma$ whose corresponding~$\qbold$ lies 
in~$\{\pm\ibold,\pm\jbold,\pm\kbold\}$.  Consequently, the quotient 
map~$\tilde N^*\to N^*$ defines a 4-fold cover of~$N^*$ that is, itself,
a 2-fold quotient of~$\hat N^*$.  I.e., there is a sequence of coverings
$$
\CD
\hat N^* @>{2{-}1}>> \tilde N^* @>{4{-}1}>>  N^*
\endCD
$$ 
corresponding to the inclusion of subgroups~$\{\bone\} \subset 
 \{\pm\bone\} \subset  \{\pm\bone,\pm\ibold,\pm\jbold,\pm\kbold\}$.  
The commutator loop~$Y^{-1}{\ast}X^{-1}{\ast}Y{\ast}X$
is closed in~$\tilde N^*$ and this is a loop over which the 
cover~$\hat N^*\to \tilde N^*$ is non-trivial since this loop does not lie 
in~$\hat K$.

It will be useful to construct a embedding of~$\tilde N^*$ 
into~$\bbC\bbP^2\times~(-\frac\pi2,\frac\pi2)$.
Consider the meromorphic solution~$\p$ on~$\bbC{\times}(-\frac\pi2,\frac\pi2)$
to the second order holomorphic differential equation with initial conditions
$$
\p_{zz}(z,\theta) = 6\,\sec^{2/3}\theta - 96\,\p(z,\theta)^2\,,
   \qquad \p(0,\theta) = r_2(\theta), \qquad \p_z(0,\theta) = 0.
$$
Of course, $\p$ is a version of the Weierstrass $\p$-function.  It satisfies
the first order differential equation
$$
\p_z(z,\theta)^2 = - 2\tan\theta + 12\,\sec^{2/3}\theta\,\p(z,\theta)
                      - 64\,\p(z,\theta)^3.
$$
Moreover~$\p(x,\theta)=f(x,\theta)$ when~$x$ is real and 
$\p(iy,\theta) = g(y,\theta)$ when~$y$ is real, as follows easily from
the Chain Rule.  Now, $\p$ is doubly
periodic and even:
$$
\p\bigl(z+2\rho_+(\theta)\bigr) = \p\bigl(z+2i\rho_-(\theta)\bigr)
= \p\bigl(-z\bigr) = \p\bigl(z\bigr).
$$
and also assumes the special values
$$
\p\bigl(i\rho_-(\theta)\bigr) = r_1(\theta),\quad
\p\bigl( 0) = r_2(\theta),\quad
\p\bigl( \rho_+(\theta)\bigr) = r_3(\theta),\quad
\p\bigl(\rho_+(\theta){+}i\rho_-(\theta)\bigr) = \infty.
$$
In fact,~$\p$ has a double pole at~$\rho_+(\theta){+}i\rho_-(\theta)$ and
no other singularities in the fundamental rectangle.  Moreover, $\p_z$
has simple zeros at~$0,\rho_+(\theta),i\rho_-(\theta)$ and a triple pole
at~$\rho_+(\theta){+}i\rho_-(\theta)$. 

Now consider, for each fixed~$\theta$ in the 
interval~$(-\frac\pi2,\frac\pi2)$, the quadratic form
$$
ds^2_\theta = \bigl(f(x,\theta)-g(y,\theta)\bigr)\,(dx^2+dy^2).
$$ 
By the earlier analysis of the vanishing locus of~$u:D\to\bbR$, 
this quadratic form defines a conformal pseudo-metric on~$\bbC$ 
that branches to order~$1$ at the points of the lattice
$$
\Lambda_\theta  = \bigl\{\ 2m\,\rho_+(\theta)+i\,2n\,\rho_-(\theta)
                      \ \vrule\ m,n\in\bbZ\bigr\}\subset\bbC
$$
and is periodic with respect to this lattice.  Since $ds^2_\theta$ is 
invariant under reflection in the $x$-axis and the $y$-axis, the 
lines~$x=m\rho_+(\theta)$  and $y=n\rho_-(\theta)$ for integer~$m$ 
and $n$ are geodesics in this metric.

The structure equations show that $ds^2_\theta$ has constant Gauss 
curvature~$K=16$, and so must be induced by pullback from the 
standard metric on the Riemann sphere with this curvature.
In particular, there is a meromorphic function~$w$ on 
the $z$@-plane so that
$$
f(x,\theta)-g(y,\theta) = \frac{|w'(z)|^2}{4\bigl(1+|w(z)|^2\bigr)^2}\,.
$$
The function~$w$ must ramify to order~$1$ at each of the points 
of~$\Lambda_\theta$ and must carry the geodesics $x=2m\rho_+(\theta)$ 
and $y=2n\rho_-(\theta)$ onto a single geodesic on the Riemann sphere. 
(Since they intersect at right angles in the $z$-plane and the intersection 
point is a branch point of~$w$ of order~$2$, the image geodesics must meet 
at an angle of~$\pi$ and hence must lie along the same geodesic on the 
sphere.)  This information is not enough to make the function~$w$ unique; 
it only determines~$w$ up to composition with an isometric rotation of 
the Riemann sphere.  However, adding the requirements that~$w(0) = 0$ and 
that $w''(0)$ be real and positive do make~$w$ unique, so this will be
assumed from now on.  

Because the geodesic segment~$t\rho_+(\theta)$ for~$0\le t\le 2$ is congruent 
to the geodesic segment~$t\rho_+(\theta)+2i\rho_-(\theta)$ for~$0\le t\le 2$,
and because the geodesic segment~$i\,t\rho_-(\theta)$ for~$0\le t\le 2$ 
is congruent to the geodesic segment~$2\rho_+(\theta)+i\,t\rho_-(\theta)$ 
for~$0\le t\le 2$, and because there are no ramification points of~$w$ in the 
interior of the fundamental rectangle, it follows that the normalized~$w$ must
satisfy
$$
w\bigl(2\rho_+(\theta)\bigr) w\bigl(2i\rho_-(\theta)\bigr) = -1
$$
with~$w\bigl(2\rho_+(\theta)\bigr)$ real and positive
and~$w\bigl(2\rho_+(\theta){+}2i\rho_-(\theta)\bigr)=\infty$.  Pursuing this
analysis, it follows without much difficulty that~$w$ must be doubly periodic
with periods~$4\rho_+(\theta)$ and $4i\rho_-(\theta)$ and have one double pole
at~$2\rho_+(\theta){+}2i\rho_-(\theta)$ in the fundamental rectangle 
of~$2\Lambda_\theta$.

By the usual properties of doubly periodic meromorphic functions
on the plane, only one function~$w$ with all these properties exists.
It can be written in terms of~$\p$ as
$$
w(z,\theta) = \frac{\p({\frac12}z,\theta)-r_2(\theta)}
{\sqrt{\bigl(r_3(\theta)-r_2(\theta)\bigr)\bigl(r_2(\theta)-r_1(\theta)\bigr)}}.
$$
The Weierstra\ss-type equation for~$\p$ shows that~$w$ itself satisfies
the Weierstra\ss-type equation
$$
(w_z)^2 = 16\,b(\theta)\,w -48\,r_2(\theta)\,w^2 -16\,b(\theta)\,w^3
$$
where 
$$
b(\theta) 
= \sqrt{\bigl(r_3(\theta)-r_2(\theta)\bigr)\bigl(r_2(\theta)-r_1(\theta)\bigr)}>0.
$$

By symmetry considerations, $w$ must map the boundary of the rectangle~$\sR$
with vertices~$0$, $2\rho_+(\theta)$, $2\rho_+(\theta){+}2i\rho_-(\theta)$,
and~$2i\rho_-(\theta)$ to the real line plus~$\infty$ on the Riemann
sphere and do so in a one-to-one and onto manner. Consequently~$w$
establishes a biholomorphism between the interior of~$\sR$ and the
upper half plane.  Because of the symmetry of the boundary values,
particularly the 
identity~$w\bigl(2\rho_+(\theta)\bigr) w\bigl(2i\rho_-(\theta)\bigr) = -1$,
it follows that~$w$ must map~$\rho_+(\theta){+}i\rho_-(\theta)$, the
center of~$\sR$, to the center of the upper half plane (endowed with its
usual metric of constant {\it positive\/} curvature), i.e., that
$$
w\bigl(\rho_+(\theta){+}i\rho_-(\theta),\theta\bigr) = i.
$$
Using this information, it is not difficult to deduce that
$$
w_z\bigl(\rho_+(\theta){+}i\rho_-(\theta),\theta\bigr) 
= 4\sqrt{3r_2(\theta) + 2b(\theta)i}.
$$
(In view of the Weierstra\ss{} equation, the only problem is to fix 
the ambiguity of the sign of this square root, but this is not
difficult.  I mean the one with positive imaginary part.)
\smallskip
It follows that there is a well-defined
map~$\Psi:\tilde N\to \bbC\bbP^2{\times}(-\frac\pi2,\frac\pi2)$ 
satisfying
$$
\Psi\bigl(\la x,y,\theta\ra\bigr) 
= \bigl([1,w(x{+}iy,\theta),w_z(x{+}iy,\theta)],\theta\bigr).
$$
Note, in particular, 
that~$\Psi\bigl(\la \rho_+(\theta),\rho_-(\theta),\theta\ra\bigr) 
= \bigl(\bigl[1,i,4\sqrt{3r_2(\theta) + 2b(\theta)i}\,\bigr],\theta\bigr)$.
The image of~$\Psi$ is the
locus~$\tilde E\subset\bbC\bbP^2{\times}(-\frac\pi2,\frac\pi2)$ consisting 
of points~$\bigl([Z_0,Z_1,Z_2],\theta\bigr)$ that satisfy the equation
$$
Z_0{Z_2}^2 = 
16\,b(\theta)\,{Z_0}^2Z_1 - 48\,r_2(\theta)\,Z_0{Z_1}^2 
       -16\,b(\theta)\,{Z_1}^3.
$$
Let~$\tilde E_\theta\subset\bbC\bbP^2$ be the smooth plane cubic curve 
so that~$\tilde E_{\theta}\times\{\theta\}
= \tilde E\cap\bigl(\bbC\bbP^2{\times}\{\theta\}\bigr)$.  This is an elliptic 
curve and will be referred to as the $\theta$@-slice of~$\tilde E$.  
By the discussion already given plus elementary properties of elliptic 
curves,~$\Psi$ is a diffeomorphism from~$\tilde N$ to~$\tilde E$.
Moreover,~$\Psi(\tilde N^*)=\tilde E^*$, which is defined
as the complement in~$\tilde E$ of the three points on each 
$\tilde E_\theta$ that lie on the line~$Z_2=0$
together with the point at infinity (i.e., the flex tangent 
on the line~$Z_0=0$) on each $\tilde E_\theta$.

Now, the double cover~$\hat N^*\to\tilde N^*\simeq \tilde E^*$ is nontrivial
around each of these missing points in each $\theta$-slice.   
Consider the smooth plane quartic 
family~$\hat E\subset \bbC\bbP^2{\times}(-\frac\pi2,\frac\pi2)$ 
consisting of points~$\bigl([W_0,W_1,W_2],\theta\bigr)$ 
that satisfy the equation
$$
{W_2}^4 = 
16\,b(\theta)\,{W_0}^3W_1 - 48\,r_2(\theta)\,{W_0}^2{W_1}^2 
       -16\,b(\theta)\,W_0{W_1}^3.
$$
The map that takes $\bigl([W_0,W_1,W_2],\theta\bigr)\in \hat E$
to~$\bigl([(W_0)^2,W_0W_1,(W_2)^2],\theta\bigr)\in \tilde E$ 
is a branched double cover over each $\tilde E_\theta$.
The branch locus over each $\tilde E_\theta$ consists of the four points 
on~$\tilde E_\theta$ that do not belong to~$\tilde E^*$.  
Let~$\hat E^*\subset\hat E$ be the inverse image of~$\tilde E^*$
under this smooth mapping.

Now the double cover~$\hat E^*\to \tilde E^*\simeq\tilde N^*$ is 
nontrivial exactly along the same curves as the double 
cover~$\hat N^*\to \tilde N^*$.  Thus, there is a 
diffeomorphism~$\hat\Psi:\hat N^*\to \hat E^*$ that identifies the
two double covers and this~$\hat\Psi$ is unique up to composition
with the deck transformation~$\bigl([W_0,W_1,W_2],\theta\bigr)\to 
\bigl([W_0,W_1,-W_2],\theta\bigr)$ 
of the covering~$\hat E^*\to\tilde E^*$.  From now on, I will
fix a choice of~$\hat\Psi$ and use it to identify~$\hat N^*$
with~$\hat E^*$. 

Each of the $\theta$-slices~$\hat E_\theta\subset \hat E$ is a nonsingular 
plane quartic and hence is a nonhyperelliptic Riemann surface of 
genus~$3$~\cite{GH, Chapter~2}.  In fact, the
functions
$$
w = \frac{W_1}{W_0}\,,\qquad v = \frac{W_2}{W_0}
$$
are smooth and well-defined on~$\hat E^*$, restricting to each~$\hat E_\theta$
to become meromorphic functions with poles located at the point at 
`infinity' given by the intersection of~$\hat E_\theta$ with the line~$W_0=0$.
The 1-forms
$$
\alpha_1 = \frac{w\,dw}{v^3}\,,\qquad 
\alpha_2 = \frac{   dw}{v^3}\,,\qquad
\alpha_3 = \frac{v\,dw}{v^3}= \frac{dw}{v^2}
$$
restrict to each~$\hat E_\theta$ to be a basis for the holomorphic
1-forms on~$\hat E_\theta$.  Note that~$\alpha_3$ is actually invariant
under the deck transformation~$(w,v,\theta)\mapsto(w,-v,\theta)$ 
of the covering~$\hat E^*\to\tilde E^*$ and hence is well-defined
as a 1-form on~$\tilde E^*$.  This 1-form restricts to each~$\tilde E_\theta$ 
to become the nontrivial holomorphic differential on that elliptic curve.  
Note that $\alpha_1$ and $\alpha_2$ have no common
zeroes:  In fact,~$\alpha_2$ has only one zero, which is of order~$4$,
and this occurs at the common pole of~$w$ and~$v$.  Since~$w$ has a pole
of order exactly~$4$ at this point, it follows that~$\alpha_1$ does
not vanish there.  

Let~$\hat n(\theta) =\bigl(\bigl[1,\,i,\,
    2 \root{4}\of{3r_2(\theta) + 2b(\theta)i}\,\bigr],\theta\bigr)$ and
consider the multivalued `function' on~$\hat E$ `defined' 
by the abelian integral
$$
\align
\vartheta([1,w,v],\theta) 
  = \pmatrix \vartheta_1([1,w,v],\theta)\\ 
             \vartheta_2([1,w,v],\theta)\endpmatrix
  &=  \int_{\hat n(\theta)}^{([1,w,v],\theta)} 
         \sqrt2\,\pmatrix w\\1\endpmatrix\frac{dw}{v^3}\\
  &=  \int_{\hat n(\theta)}^{([1,w,v],\theta)} 
         \pmatrix \sqrt2\,\alpha_1\\ \sqrt2\,\alpha_2 \endpmatrix.\\
\endalign
$$
where the integral is to be computed along a path joining $\hat n(\theta)$
to~$([1,w,v],\theta)\in \hat E$ that lies entirely in~$\hat E_\theta$.
Of course, the value of this integral depends on the homology class
of the path joining the two endpoints, so this is not well-defined
as a function on~$\hat E$.  The ambiguity in the definition of~$\vartheta$ 
will be determined below.  For the time being, consider~$\vartheta$ as being 
defined on a suitable cover~$\check E\to\hat E$.  Since~$\alpha_1$
and $\alpha_2$ do not have any common zeroes, this map is an immersion
on each $\hat E_\theta$.

Now consider the functions
$$
A = \frac{\bar v}{v\sqrt{1+|w|^2}},\qquad\qquad
B = \frac{-\bar v\,w}{v\sqrt{1+|w|^2}}
$$
defined on~$\hat E^*$.  They satisfy~$|A|^2+|B|^2 =1 $, so the function
$$
h = \pmatrix 
1 & 0 & 0 \\ 
0 & \bar A\bigl(\hat n(\theta)\bigr) & - B\bigl(\hat n(\theta)\bigr)\\ 
0 & \bar B\bigl(\hat n(\theta)\bigr) &   A\bigl(\hat n(\theta)\bigr)\\ 
\endpmatrix^{-1}
\pmatrix 1 & 0 & 0 \\ 
             \vartheta_1 & \bar A & - B\\ 
             \vartheta_2 & \bar B & A\\ \endpmatrix
$$
takes values in~$G_0$ and
is well-defined on the open set~$\check E^*\subset \check E$ that is
the inverse image of~$\hat E^*$ under the cover~$\check E \to \hat E$.  
(The purpose of the first matrix is to
arrange that~$h\bigl(\hat n(\theta)\bigr) = \text{I}_3$ for all~$\theta$,
which will be needed below.)

Since the first factor in~$h$ depends only on~$\theta$, computation yields
$$
h^{-1}\,dh \equiv
\pmatrix 
0 & 0 & 0\\
A\,d\vartheta_1+B\,d\vartheta_2 & A\,d\bar A + B\,d\bar B & B\,dA - A\,dB\\
\bar A\,d\vartheta_2-\bar B\,d\vartheta_1 
           & \bar A\,d\bar B - \bar A\,d\bar A & \bar A\,dA + \bar B\,dB\\
\endpmatrix \bmod d\theta.
$$
Since $d\vartheta_i\equiv \sqrt2 \alpha_i\bmod d\theta$ 
for~$i=1$, $2$, it follows that
$$
\left.
\aligned
A\,d\vartheta_1+B\,d\vartheta_2 &\equiv 0 \\
\bar A\,d\vartheta_2-\bar B\,d\vartheta_1 &\equiv \omega \\
\endaligned
\right\}\quad \bmod d\theta.
$$
(This last follows from the identities~
$$
\omega \equiv (f-g)^{-1/4}\,dz 
\equiv \frac{\sqrt{2(1+|w|^2)}}{|w'(z)|^{1/2}} \frac{dw}{w'(z)}
\equiv \frac{\sqrt{2(1+|w|^2)}}{|v|} \frac{dw}{v^2} \quad\bmod d\theta,
$$
together with the definitions of~$A$ and $B$.  The reader can now probably
see why the factor of~$\sqrt2$ was introduced into the definition 
of~$\vartheta$.)  Moreover, 
$$
\align
\bar A\,d\bar B - \bar B\,d\bar A 
&= -{\bar A}^2\,d(\bar B/\bar A) 
=  \frac{-v^2}{|v|^2(1+|w|^2)}\,d\bar w\\
&= \frac{-|v|^2}{(1+|w|^2)}\,\frac{d\bar w}{{\bar v}^2}
 \equiv -2 s\,\omega \equiv \sigma \mod d\theta.
\endalign
$$
By these results, there exists a real-valued
1-form~$\phi^*$ so that
$$
h^{-1}\,dh \equiv 
\pmatrix 0 & 0 & 0\\
         0 & i\,\phi^* & -\bar\sigma\\
        \omega & \sigma & - i\,\phi^*\\
\endpmatrix \bmod d\theta.
$$
The matrix on the right is almost $\gamma$.  In fact, I claim that
it is congruent to~$\gamma$ modulo~$d\theta$.  Since~$\eta$ is a multiple
of~$d\theta$ by definition, the only thing to check is whether~$\phi^*\equiv
\phi$ modulo~$d\theta$.  However, this follows immediately from the 
structure equations, which show that~$d\sigma \equiv 2i\,\phi\w\sigma \bmod
d\theta$ while the very fact that $\phi^*$ appears where it does 
in~$h^{-1}\,dh$ shows that $d\sigma \equiv 2i\,\phi^*\w\sigma \bmod
d\theta$.  Comparing these two relations and using the fact that~$\phi$ and
$\phi^*$ are real then yields~$\phi^*\equiv\phi\bmod d\theta$, as desired.
(Alternatively, one can simply carry out the computations and compare
the results.)

It has now been shown that $h^{-1}\,dh \equiv \gamma \mod d\theta$.
Now, $\gamma$ is well-defined on~$\hat E^*$, not just on~$\check E^*$,
so it follows that~$h^{-1}\,dh$ is well defined on each~$\hat E_\theta$
and has the same holonomy as~$\gamma$ on each~$\hat E_\theta$. 
Now, it has already been shown that the holonomy of~$\gamma$ on~$\hat E^*$ 
lies in the discrete subgroup~$\hat\Lambda\subset G_0$ and the inclusion
$\hat E^*_\theta\hookrightarrow \hat E^*$ induces and isomorphism on 
fundamental groups. Consequently, there is a well-defined mapping
$$
\hat\Lambda h: \hat E \to \hat\Lambda\backslash G_0.
$$ 
Note that the quotient is via the left action and not the right action.  
In particular, the canonical left-invariant form on~$G_0$ is well-defined 
on $\hat\Lambda\backslash G_0$.

Now, consider the $\eug_0$-valued 1-form~$\kappa$ that is 
well-defined on~$\check E^*$ by the formula
$\kappa = h\,\gamma\,h^{-1} - dh\,h^{-1}$.
Since
$$
    \kappa   = h\,(\gamma-h^{-1}\,dh)\,h^{-1}\equiv0 \bmod d\theta
$$
since~$d\kappa = -\kappa\w\kappa = 0$, and since~$\kappa$ vanishes
when restricted to each~$\check E^*_\theta$, it must be a 1-form in~$\theta$
alone.  In fact, a computation using the properties of~$f$ and~$g$ shows
that
$$
\hat n^*(\kappa) = \pmatrix 0 & 0 & 0 \\ 
\frac16(r_3(\theta)-r_1(\theta))^{1/4}\sec^{2/3}\theta\,d\theta & 0&0\\
0&0&0\\
\endpmatrix.
$$
In particular, $\kappa= k^{-1}\,dk$ where
$$
k = \pmatrix 1 & 0 & 0\\
            m(\theta) & 1 & 0\\
             0 & 0 & 1\\
     \endpmatrix.
$$
and where~$m$ satisfies~$m(0)=0$ and $m'(\theta) 
= \frac16(r_3(\theta)-r_1(\theta))^{1/4}\sec^{2/3}\theta$.  

Since the elements of the form~$k(\theta)$ commute with all of the
elements of~$\hat\Lambda$, it now follows that $\gamma = g^{-1}\,dg$
where $[g] = \hat\Lambda kh$ is well defined on~$\hat E^*$ as a map 
into~$\hat\Lambda\backslash G_0$.  Since~$\hat\Lambda\backslash G_0/K\simeq
\bbC^2/\Lambda$, the map
$$
\Phi([1,w,v],\theta) \equiv 
\pmatrix 
\phantom{-}A\bigl(\hat n(\theta)\bigr) &  B\bigl(\hat n(\theta)\bigr)\\ 
-\bar B\bigl(\hat n(\theta)\bigr) &  \bar A\bigl(\hat n(\theta)\bigr)\\ 
\endpmatrix
\pmatrix \vartheta_1([1,w,v],\theta) \\ \vartheta_2([1,w,v],\theta)
\endpmatrix
+ \pmatrix m(\theta) \\ 0 \endpmatrix \mod\Lambda
$$
is well-defined as a map~$\Phi:\hat E \to \bbC^2/\Lambda$.

From the formulae that went into its definition, $\Phi$ is an immerison 
on~$\hat E^*$ whose image is a Levi-flat minimal hypersurface 
in~$\bbC^2/\Lambda$ of type~$(3)$.  Moreover~$\Phi\bigl(\hat E_\theta\bigr)
\subset\bbC^2/\Lambda$ is a complex leaf in this hypersurface and
is immersed as a compact Riemann surface of genus~$3$.

Now, $\Phi$ is not an immersion near the four curves~$v=0$ in~$\hat E$. 
(These are the curves that intersect each~$\hat E_\theta$ in the four
branch points.)  In fact, it collapses each of these curves to a point,
as can be seen by doing a local computation.  Let these points be labeled
$P_i\in\bbC^2/\Lambda$ for~$i=1$ to $4$.   

\remark{A possible `algebraic' structure}
Now~$\bbC^2/\Lambda$ is a complex torus that has nontrivial divisors, 
for example, the genus~$3$ Riemann surfaces~$\Phi\bigl(\hat E_\theta\bigr)$.
It follows that~$\bbC^2/\Lambda$ is an Abelian variety (actually, this
also follows from the explicit description of~$\Lambda$ as a lattice of
type~$F_4$ that has already been given).  In particular, $\bbC^2/\Lambda$ 
is an algebraic surface. By a standard Riemann-Roch 
calculation~\cite{GH, Chapter 4}, one can show that the curves in the 
connected family of~$C_\theta = \Phi\bigl(\hat E_\theta\bigr)$ that pass
through the points~$P_i$ form a pencil, i.e., the moduli~$M$ of such 
curves is a~$\bbC\bbP^1$.  In fact, regarding~$\theta$ as
a complex parameter in the formula for~$\Phi$ gives a local real
parameter on~$M$ near~$\theta=0$.  Evidently, the curve~$M$ admits 
an antiholomorphic involution for which 
the curves~$C_\theta$ are fixed points.  Of course, any antiholomorphic
involution of~$\bbC\bbP^1$ that has fixed points is conjugate via an
automorphism of~$\bbC\bbP^1$ to the standard conjugation fixing 
an~$\bbR\bbP^1\subset\bbC\bbP^1$.  Thus, it would appear that the 
image~$\Sigma=\Phi\bigl(\hat E\bigr)\subset \bbC^2/\Lambda$ is a dense
open set in an `algebraic' real hypersurface~$\bar\Sigma\subset
\bbC^2/\Lambda$ that is the union of the curves in~$M$ that are fixed
under the antiholomorphic involution.
Presumably, the singular curves in the pencil~$M$ are
unions of elliptic curves embedded in~$\bbC^2/\Lambda$ linearly and
are therefore the totally geodesic complex leaves in~$\bar\Sigma$. It would 
be interesting to know whether or not the only singularities of~$\bar\Sigma$ 
are the four points~$P_i$ and whether or not these singular points really 
do resemble cones on the Clifford torus, as they appear to. 
\endremark

\enddocument